\documentclass[lfeqn,12pt]{article}
\usepackage{amssymb,bbm,pifont,mathtools,mathabx,mathrsfs, amsthm, amsmath%fdsymbol,
}

\usepackage{times}
\usepackage[nottoc]{tocbibind}

\usepackage{setspace} 
 
\usepackage{comment}
\usepackage{multirow}
\usepackage[multiple]{footmisc}
\usepackage{xcolor}
\usepackage[colorlinks=true]{hyperref}
\usepackage{marginnote}
\newtheorem{theorem}{Theorem}[section]
\newtheorem{definition}[theorem]{Definition}
\newtheorem{proposition}[theorem]{Proposition}
\newtheorem{lemma}[theorem]{Lemma}
\newtheorem{algorithm}[theorem]{Algorithm}
\newtheorem{remark}[theorem]{Remark}
\newtheorem{sublemma}[theorem]{Sublemma}
\newtheorem{corollary}[theorem]{Corollary}
\newtheorem{problem}[theorem]{Problem}
\newtheorem{example}[theorem]{Example}
\newtheorem{assumption}[theorem]{Assumption}
\newtheorem{notationalrem}[theorem]{Notational Remark}
\newtheorem{tools}[subsection]{\( \negsp\negsp \)}
\newcommand\asm[1]{ \begin{assumption}\label{#1} }
\newcommand\easm{ \end{assumption} }
\newcommand\dfn[1]{ \begin{definition}\label{#1} }
\newcommand\dfntwo[2]{ \begin{definition}[#2]\label{#1} }
\newcommand\edfn{ \end{definition} }
\newcommand\rem[1]{ \begin{remark}\label{#1} \small \rm}
\newcommand\remtwo[2]{ \begin{remark}[#2]\label{#1} \rm}
\newcommand\erem{ \end{remark} }
\newcommand\thm[1]{ \begin{theorem}\label{#1}}
\newcommand\thmtwo[2]{ \begin{theorem}[#2]\label{#1}}
\newcommand\ethm{ \end{theorem} }
\newcommand\pro[1]{ \begin{proposition}\label{#1}}
\newcommand\protwo[2]{ \begin{proposition}[#2]\label{#1}}
\newcommand\epro{ \end{proposition} }
\newcommand\lem[1]{ \begin{lemma}\label{#1}}
\newcommand\lemtwo[2]{ \begin{lemma}[#2]\label{#1}}
\newcommand\elem{ \end{lemma} }
\newcommand\sublem[1]{ \begin{sublemma}\label{#1}}
\newcommand\sublemtwo[2]{ \begin{sublemma}[#2]\label{#1}}
\newcommand\esublem{ \end{sublemma} }
\newcommand\cor[1]{ \begin{corollary}\label{#1}}
\newcommand\cortwo[2]{ \begin{corollary}[#2]\label{#1}}
\newcommand\ecor{ \end{corollary} }
\newcommand\notrem[1]{ \begin{notationalrem}\label{#1} \sl}
\newcommand\enotrem{ \end{notationalrem} }

            \newcommand\prob[1]{ \begin{problem}\label{#1} }
\newcommand\eprob{ \end{problem} }
\newcommand\examp[1]{ \begin{example}\label{#1} }
\newcommand\eexamp{ \end{example} }
\newcommand\beq[1]{ \begin{equation}\label{#1} }
\newcommand{\eeq}{ \end{equation} }

\newcommand\beqa[1]{ \begin{eqnarray} \label{#1}}
\newcommand{\eeqa}{ \end{eqnarray} }
\newcommand{\beqano}{ \begin{eqnarray*} }
\newcommand{\eeqano}{ \end{eqnarray*} }

\newcommand\ovl[1]{ \overline {#1} }

\newcommand\inter{ {\, \rm int\, }}

\renewcommand{\=}{ {\, \equiv\, } }

\newcommand{\negsp}{\hspace{-.04truecm}}
\font\teneufm=eufm10
\font\seveneufm=eufm7
\font\fiveeufm=eufm5
\newfam\eufmfam
\textfont\eufmfam=\teneufm
\scriptfont\eufmfam=\seveneufm
\scriptscriptfont\eufmfam=\fiveeufm

\newcommand{\kw}[1]{\textbf{\emph{#1}}}
\usepackage{authblk}
\begin{document}
%%%%%%%%%%%%%%Header
%%%%%%%%%%%%%%%
\date{\today}

\title{Rigorous numerics for critical orbits \\ in the quadratic family}

\author{A. Golmakani\footnote{Ali Golmakani, Universidade Federal de Alagoas, Av.\ Lourival Melo Mota, s/n, Macei\'{o}, Alagoas 57072-900, Brazil;
\emph{aligolmakani@gmail.com}}, 
\
C. E. Koudjinan\footnote{Comlan Edmond Koudjinan, Abdus Salam
International Centre for Theoretical Physics (ICTP),
Strada Costiera 11, 34151 Trieste, Italy; \emph{koudjinanedmond@gmail.com}
}, 
\
S. Luzzatto \footnote{Stefano Luzzatto, Abdus Salam
International Centre for Theoretical Physics (ICTP),
Strada Costiera 11, 34151 Trieste, Italy; \emph{luzzatto@ictp.it}
},
\ 
P. Pilarczyk\footnote{Pawe\l{} Pilarczyk, Gda\'{n}sk University of Technology,
Faculty of Applied Physics and Mathematics,
ul.\ Gabriela Narutowicza 11/12, 80-233 Gda\'{n}sk, Poland;
\emph{pawel.pilarczyk@pg.edu.pl}}
}
\date{}
\maketitle

\begin{abstract}
We develop algorithms and techniques to compute \emph{rigorous} bounds for finite pieces of orbits of the critical points, for \emph{intervals of parameter values}, in the quadratic family of one-dimensional maps \(f_a (x) = a - x^2\). We illustrate the effectiveness of our approach by constructing a dynamically defined partition \( \mathcal P \) of the parameter interval \( \Omega=[1.4, 2] \) into almost 4 million subintervals, for each of which we compute to high precision the orbits of the critical points up to some time~\( N \) and other dynamically relevant quantities, several of which can vary greatly, possibly spanning several orders of magnitude. We also subdivide \( \mathcal P \) into a family \( \mathcal P^{+} \) of intervals which we call \emph{stochastic intervals} and a family \( \mathcal P^{-} \) of intervals which we call \emph{regular intervals}. We numerically prove that each interval \( \omega \in \mathcal P^{+} \) has an \emph{escape time}, which roughly means that some iterate of the critical point taken over all the parameters in \( \omega \) has considerable width in the phase space. This suggests, in turn, that most parameters belonging to the intervals in \( \mathcal P^{+} \) are stochastic and most parameters belonging to the intervals in \( \mathcal P^{-} \) are regular, thus the names. We prove that the intervals in \( \mathcal P^{+} \) occupy almost 90\% of the total measure of \( \Omega \). 
The software and the data is freely available at \href{http://www.pawelpilarczyk.com/quadr/}{http://www.pawelpilarczyk.com/quadr/},
and a web page is provided for carrying out the calculations.
The ideas and procedures can be easily generalized to apply to other parametrized families of dynamical systems.
\end{abstract}

\vskip 24pt

\noindent
In the 1970s Robert May introduced the logistic family of one-dimensional maps as an example of a simple mathematical model which nevertheless exhibits extremely complex behaviour. Since then, the logistic family and the very closely related quadratic family have become an icon of Chaos Theory. Notwithstanding some very deep analytic and abstract results obtained over the last several decades by top mathematicians, and extensive numerical studies by physicists and nonlinear dynamicists, starting from Feigenbaum, there are literally only a handful of rigorous concrete numerical results. This is not too surprising because it is indeed the essence of the chaotic dynamics on these families which makes them numerically very challenging.

In this paper we develop some rigorous numerical techniques for studying the quadratic family and obtain several interesting ``statistical'' results about how often certain dynamical situations occur in parameter space. In particular, we conclude that stochastic-like dynamics is likely to occur for almost 90\% of parameters. Our research is motivated by a specific ambitious project to identify true chaotic dynamics in the family. However, our techniques can certainly be easily adapted to a large variety of situations.

%{\bf Keywords:}

%\setcounter{tocdepth}{1}
%\tableofcontents

% ====================================================

\section{Introduction}
\label{sec:intro}

The rigorous computation of orbits of dynamical systems is well known to be very delicate due to inevitable approximation errors caused by the fact that computers work with only a finite set of ``representable'' numbers, such as the 64-bit floating point numbers following the IEEE 754 standard \cite{ieee754}, implemented in most modern processors. A standard and effective way to deal with this problem is to use \emph{interval arithmetic}~\cite{Moore1966,WT2011} to obtain rigorous bounds for the iterates of a single point, which can be made arbitrarily sharp by paying the price in computing time. The situation can, however, get significantly more complicated if we need to bound the images of an ``ensemble'' of points or the images of a single point for different parameter values. The purpose of this paper is to illustrate some of the problems and provide computational techniques to address them. We focus on a particular case which is motivated by a bigger and more ambitious project, as explained below. However, similar problems appear in more general situations, and our approach should be relatively straightforward to apply in other settings.

% ====================================================

\subsection{The quadratic family}

We consider the classical quadratic family of one dimensional maps given by 
\begin{equation}\label{eq:quadratic}
f_{a}(x) = a-x^{2}
\end{equation}
and restrict ourselves to parameters 
\( a\in \Omega:=[1.4, 2] \), since the dynamics of \( f_{a} \) is essentially trivial and well understood for \( a\notin\Omega \),
and initial conditions \( x\in I_{a} \), 
where the interval~\( I_{a} \) depends continuously on the parameter \(  a  \) and has the property that \(  f(I_{a})\subseteq I_{a}  \), and that the iterates of all the points \(  x\notin I_{a}  \) converge to \(  -\infty  \). The existence of \( I_{a} \) follows by elementary observations and its properties imply that  any non-trivial dynamics is contained in \(  I_{a}  \).

Let \(  \omega\subseteq \Omega  \) be an arbitrary parameter interval. Formally, we could even take \(  \omega=\Omega  \), but in general our calculations are most effective for quite small intervals. In the computations to be given below as an illustration of our methods, we will construct dynamically a partition of \( \Omega \) into subintervals \( \omega \) whose length varies from an order of \( 10^{-3} \) to as small as \( 10^{-10} \).  
Let \( c \) denote the critical point \( 0 \) of \( f_a \). For each \(  n\geq 0  \), we let 
\begin{equation}\label{eq:omegan}
c_{n}(a):= f_{a}^{n}(f_{a}(c)) \qquad \text{ and } \qquad 
\omega_{n}:=\{c_{n}(a): a\in \omega\}.
\end{equation}
Notice that the critical value \(  c_{0}(a) \) equals \( a  \); therefore, \(  \omega_{0}  \) coincides with \( \omega  \). For \(  n\geq 1  \), \(  c_{n}(a)  \) is simply the \(  n  \)'th image of the critical value and \(  \omega_{n}  \) is the interval given by the \(  n  \)'th images of the critical values for all the parameters \(  a\in \omega  \). 

The first and main objective of this paper is to describe and implement effective computational techniques to obtain \emph{arbitrarily sharp} and \emph{rigorous} approximations  for  \(  \omega_{n}  \) under a verifiable technical assumption \eqref{eq:mon} to be given below. We will also describe arguments to obtain rigorous bounds on a few other relevant dynamical quantitities. These objectives are motivated by a bigger project that we discuss in the following subsections. In Section \ref{sec:results} we present and discuss the results arising from our computations. 
In Section \ref{sec:strategies} we give a relatively detailed overview of the computational strategies used to achieve our goals, and in Section \ref{sec:P} explain how these are used to construct the dynamically defined partition \( \mathcal P \). In Section~\ref{sec:procedure} we give all the details of the computational procedures  and explain how we are able to ensure rigorous bounds, and in  Section~\ref{sec:algorithms} we give details of the algorithms.
The source code of the software, programmed in C++, is freely available at the website \cite{software}, which also features a user-friendly interface to run the software directly from the web browser.
The data resulting from our computations is published in~\cite{data20}.

% ====================================================

\subsection{Regular and stochastic dynamics}
\label{sec:motivation}
The specific approach developed in this paper concerns calculations of quantities of very general interest, in a variety of settings relevant to anyone studying dynamical systems from a numerical point of view. In our case, they are directly motivated by a more ambitious long-term research programme whose main interest lies precisely in the subtle and non-trivial synergy between rigorous computational methods and more standard analytic, geometric and probabilistic mathematical arguments. In this section we outline the main features and goals of this programme and emphasize the crucial role of the computational methods introduced in this paper. 

The quadratic family \eqref{eq:quadratic} of one-dimensional maps is possibly one of the most studied families of dynamical systems. It contains a mind-boggling richness of dynamical phenomena, which has still not been completely classified or understood,  and the dependence of the dynamics on the parameter  is extremely complicated. It is known, however, that only two types of dynamical phenomena occur with  positive probability in the parameter interval \( \Omega \): \emph{regular} dynamics, where \(  f_{a}  \) admits a unique attracting periodic orbit to which Lebesgue almost every \(  x\in I_{a}  \) converges, or \emph{stochastic} dynamics, where \(  f_{a}  \) admits a unique invariant probability measure \(  \mu_{a}  \) to which the ergodic averages of Lebesgue almost every point \(  x\in I_{a}  \) converge (in a very ``chaotic'' way, thus the term ``stochastic-like''). In other words, the union of the two sets 
 \begin{equation}\label{eq:regstoc}
\Omega^{-}:=\{a\in \Omega: a \text{ regular}\}
\quad \text{ and } \quad 
\Omega^{+}:=\{a\in \Omega: a \text{ stochastic}\},
\end{equation}
has  \emph{full measure} in \( \Omega \)  \cite{Lyu02, AviLyudMel03}.  It is also known  that \(  \Omega^{-}  \) is \emph{open and dense} in \(  \Omega  \) \cite{GraSwi97, Lyu97, Lyu97a} and therefore  \(  \Omega^{+}  \)  is \emph{nowhere dense}, but has \emph{positive Lebesgue measure} \cite{Jak81, BenCar85}. 
A natural question is:
\begin{quote}
Given an explicit parameter \( a\in \Omega^{-}\cup\Omega^{+} \), can we decide if \( a\in \Omega^{-} \) or \( a\in \Omega^{+} \)?
\end{quote}
It turns out that for most parameters in  \( \Omega^{-}\cup\Omega^{+} \), this is an extremely difficult question, and the set  \( \Omega^{+} \) is  in fact formally \emph{undecidable} \cite{ArbMat04}.   
 Nevertheless, some results do exist. Rigorous computer assisted arguments have been developed in \cite{TucWil09} to explicitly compute intervals of parameters belonging to \( \Omega^{-} \). These arguments have been applied, at the cost of an equivalent of a whole year of CPU time,  to the logistic family \( g_{\lambda}(x)=\lambda x (1-x) \) to show that at least \( 10.2\% \) of parameters in a parameter interval roughly corresponding to our interval \( \Omega \) belong to \( \Omega^{-} \); these parameters apparently consist of almost 5~million subintervals corresponding to regions with associated attracting periodic orbits of period up to about \( 30{,}000 \). An improved method was later applied in \cite{Gal17} to obtain a slightly better estimate with considerably lower computation time. The logistic family  is in fact smoothly conjugate by an explicit formula to the quadratic family \eqref{eq:quadratic} and so in principle the periodic windows for the quadratic family can be known explicitly by taking images of those computed for the logistic family.  Since the conjugacy is nonlinear, an estimate of the corresponding measure is non trivial, and will be computed in a future paper, though it turns out to yield very similar estimates, thus leaving almost \( 90\% \) of parameters unaccounted for; indeed, the results we present below are very much aligned with this figure.

Approaching the problem from the other side, notwithstanding the impossibility in general to establish that a given parameter \( a\) belongs to  \(\Omega^{+} \), it may be possible to assign a well-defined lower bound to the  \emph{probability} that \( a\in \Omega^+ \). Suppose, for example, that \( \omega \) is a small neighbourhood of the parameter \( a \) in \( \Omega \) and that, letting \( \omega^{+}:=\omega\cap\Omega^{+} \), we could show that  \( |\omega^{+}|\geq  \eta |\omega| \) for some \( \eta\in (0,1) \). Then we could say that the probability that \( a\in \omega^{+} \) is at least \( \eta \). 
The very first proof that \( |\Omega^{+}|>0 \) goes back to Jakobson \cite{Jak81}, after which there have been many generalizations  \cite{BenCar85, BenCar91,  LuzTuc99, LuzVia00, PacRovVia98}, all based on a combination of analytic, combinatorial and probabilistic arguments which imply that for some \emph{sufficiently small}  neighbourhood \( \omega \) of some ``good'' parameter value \( a^{*} \) we have \( |\omega^{+}|>0 \). However, none of the papers cited provides any explicit lower bound for the measure of \( \omega^{+} \). 

In \cite{Jak01}, Jakobson extended the arguments developed in \cite{Jak81} towards a more explicit and quantitative formulation, and designed an algorithm to estimate rigorously from below the measure of stochastic parameter values in quadratic and similar smooth families of unimodal maps. It is worth noting that this is not simply a matter of ``keeping track of the constants'' but requires a reformulation of some of the starting conditions of the results in order to make them computationally verifiable, and a corresponding modification of the arguments. The actual implementation of such an algorithm was however first carried out in \cite{LuzTak06}, using arguments more closely related to \cite{BenCar85, BenCar91}, where it was shown that  \( 97\% \) of parameters in the interval  \( \omega:= [2-10^{-4990}, 2]  \) are stochastic, thus implying that  \(  |\omega^{+}|\geq 0.97 \cdot 10^{-4990} \geq 10^{-5000}  \). This is of course an extremely small lower bound and undoubtedly very far from optimal in terms of the overall measure of stochastic parameters in  \( \Omega^{+} \) , but notwithstanding several preliminary announcements, it still remains to this day the only explicit and rigorous bound available. 
In Section \ref{sec:compstart} we briefly outline a possible strategy for extending the arguments of \cite{LuzTak06} to other parameter intervals in \( \Omega \) and explain how the results and calculations presented in this paper form part of this strategy.

% ====================================================

\subsection{Computable starting conditions}
\label{sec:compstart}
Extending the methods introduced in \cite{LuzTak06} to other parameter intervals in \( \Omega \) requires non-trivial computer-assisted calculations in order to verify some explicit \emph{starting conditions}, which were verified analytically in \cite{LuzTak06} by choosing a very small neighbourhood of the special parameter value \( a^{*}=2 \). 
It is beyond the scope of this paper to give a complete and precise list of the quantities which need to be calculated, so we refer the reader to \cite{LuzTak06} for the full technical details. Here we limit ourselves to a heuristic (and incomplete) overview which we hope nevertheless helps to get a preliminary idea and to motivate the results presented in this paper. 

We suppose first of all that we have fixed  a parameter interval \( \omega\subseteq \Omega \). Some  conditions, labelled as (A1)-(A4) and  involving a number of constants, are formulated in \cite{LuzTak06} where  it is proved that if these conditions are satisfied for a set of constants which satisfy certain inequalities, then an explicit formula gives a rigorous lower bound for the proportion of stochastic parameters in \( \omega \). A crucial and non-trivial aspect of the result is that the required conditions (A1)-(A4) are all  \emph{verifiable} and the corresponding constants are \emph{computable},  albeit by highly non-trivial computations, in finite time and with finite precision (unlike the starting conditions of the generalizations of Jakobson's Theorem mentioned above, apart from some very exceptional cases).  

The first two conditions, (A1) and (A2), are by far the most important, while (A3) and (A4) can be considered ``technical'' and may possibly even be relaxed to some extent. We therefore  focus on the first two. Without going into the precise formulation of condition (A1) we mention that it involves the choice of a  constant \( \delta>0\) which defines the critical neighbourhood 
\begin{equation}
\label{eq:Delta}
\Delta:=(-\delta, \delta). 
\end{equation}
Notice that the critical point \( c=0 \) is a critical point for all parameter values \( a \) and thus \( \Delta \) can be chosen independently of the parameter \( a \).  
Condition (A1) then essentially says that there exists a constant \( \lambda>0 \) such that the derivative \( |(f^{k})'(x)| \) of any initial condition \( x \) is \emph{growing exponentially} with exponential rate \( \lambda \), i.e. \( |(f^{k})'(x)| \geq Ce^{\lambda k} \) for some constant \( C>0 \) independent of \( x \), as long as the images of \( x \) stay outside the critical neighbourhood, i.e. as long as \( x, f(x),...,f^{k-1}(x)\notin\Delta \). This is a highly non-trivial condition if \( \delta \) is small (which it needs to be in order for the overall argument to work) since the orbit of \( x \) can still pick up some very small derivatives even outside \( \Delta \). It can be verified analytically in ``sufficiently small'' parameter neighbourhoods \( \omega  \) (whose size is however not explicitly known) of ``good'' parameters \( a^{*} \) defined by conditions which are in general also not explicitly verifiable.  The only option to verify this condition in general parameter intervals \( \omega \) is therefore by direct and explicit computation. Rigorous algorithms and computational techniques  for this purpose were developed in  \cite{DKLMOP08,GolLuzPil16} based on the construction of some relevant weighted directed graphs. 

The exponential growth of the derivative outside the critical neighbourhood \( \Delta \) is an open condition in parameter space and is in itself compatible with pretty much any kind of overall asymptotic dynamical behaviour. Indeed, as mentioned above, the set \( \Omega^{-} \) is open and dense in \( \Omega \) and therefore any interval \( \omega \) will contain a non-empty (in fact open and dense) subset of regular parameters which admit an attracting periodic orbit.  Our objective however is to show that \( \omega \) also contains stochastic parameters and indeed to obtain a lower bound for the proportion of stochastic parameters in \( \omega \). By standard results, a sufficient condition for a parameter \( a \) to be stochastic is the \emph{Collet-Eckmann} condition that the derivative along the orbit of the critical value \( c_{0}:=f_{a}(c) \) is growing exponentially fast, i.e. that there exist constants \( C, \lambda >0 \) such that \( |(f^{n})'(c_{0})| \geq Ce^{\lambda n} \) for \emph{every} \( n\geq 1 \). If condition (A1) discussed above holds, then this is satisfied as long as the orbit of the critical value stays outside \( \Delta  \) for all iterates, which can and does indeed happen but only for an exceptional set of parameters of zero Lebesgue measure. To obtain meaningful results we cannot therefore avoid having to deal with returns of the critical value to the critical neighbourhood, and in fact to returns which may come arbitrarily close to the critical point. In these cases it is impossible to verify the Collet-Eckmann condition computationally because it is not implied by any finite time condition and therefore we would need to check directly the derivative for an infinite number of iterates. We remark that there exist also weaker sufficient conditions for the parameter \( a \) to be stochastic, in some cases it is for exmaple sufficient to show that \( |(f^{n})'(c_{0})| \to \infty\), but they are still all not computationally verifiable since they are all asymptotic conditions that cannot be checked in any finite number of iterations. This is essentially the reason why stochastic parameters are undecidable, as mentioned above, and why they occur as Cantor sets and not open sets of parameters. 

The strategy, first developed by Jakobson, and refined in subsequent papers to deal with the situation described above, is to set up a \emph{probabilistic} argument based on two fundamental facts. The first, which is relatively elementary, is that the exponential growth of the derivative for the critical orbit is implied by a \emph{bounded recurrence} condition on the critical orbit, essentially something of the form \( |c_{n}|\geq e^{-\alpha n} \) for all \( n\geq 1 \) and for some sufficiently small \( \alpha>0 \) (in fact a little bit more is needed but this gives the main idea). This condition allows the critical point to be recurrent, i.e. to have arbitrarily close returns, but in a sufficiently controlled way, and also suggests that one way to establish abundance of stochastic parameters is to show that many of them have bounded recurrence. Based on this observation, the second, and much more sophisticated, key part of the strategy is to show that the intervals \( \omega_{n} \), which are precisely the union of images \( c_{n}(a) \) of the critical points for the parameters in \( \omega \), tend to \emph{grow} (exponentially fast), implying that the points \( c_{n}(a) \) are sufficiently ``spread out'' in the phase space and thus only a very small proportion can actually come close to the critical point and fail the bounded recurrence condition. 

The growth in size of the intervals \( \omega_{n} \) is thus an essential ingredient in all the proofs of all variations of Jakobson's Theorem. The proof of this fact is very involved and requires a combination of several techniques, including some combinatorial, analytic and probabilistic arguments, which themselves however  rely on features of the dynamics corresponding to the parameters in \( \omega \). It turns out that the uniform expansivity outside the critical neighbourhood \( \Delta \), as formulated in condition (A1) and as mentioned above, is one of the two most crucial features required. The second is formulated in  condition (A2) which uses the definition of \emph{escape time} which we formulate  here in a slightly simplified form as follows.
\begin{definition} \label{def:escape}
\( N \) is called an \emph{escape time} for \( \omega \) if the following holds:
\begin{equation}\label{eq:escape}
 \omega_{i}\cap \Delta=\emptyset \quad \text{for all \( 0\leq i < N \)}, \quad \text{and}  \quad |\omega_{N}|\geq \sqrt\delta. 
\end{equation}
\end{definition} 
This says  that all intervals \( \omega_{n} \) remain outside the critical neighbourhood (and thus in particular ``benefit'' from the expansivity provided by (A1)) up to time \( N \) and that they grow to ``large scale'' (in this case defined as \( \sqrt\delta \) but this can be flexible) at time~\( N \). 
The wording ``escape time'' is purposefully borrowed from  \cite{BenCar91} and later generalizations such as \cite{LuzTuc99, LuzVia00, LuzTak06}, and attempts the capture the idea, mentioned above, that the large size of the interval \( \omega_{N} \) implies that most images do not fall close to the critical point and therefore ``escape'' the constraints of the bounded recurrence condition. 

We remark that the foreseen future applications of our estimates to the general problem of the measure of stochastic parameters, and the actual formulation of condition (A2)  requires \( N \) to be ``sufficiently large'' depending  on the other constants involved, such as the size of the critical neighbourhood \( \Delta \) and the expansivity exponent \( \lambda \). 
The  main goal of this paper is for the moment more limited and is to develop the computational techniques to construct a large number of (small) parameter intervals which have an escape time at some  (possibly large) value of \( N \).
The verification of an escape time requires the computation of rigorous  \emph{enclosures} (we give the precise definitions below) of the sequence of intervals \( \omega_{i} \) for \( 0\leq i \leq N \). For these reasons, the main technical part of this paper consists of the  development of  some very efficient and effective procedures for estimating the precise location and size of intervals  \( \omega_{i} \). 

In view of future applications to parameter exclusion arguments, but also out of  independent interest, we compute some additional quantities related to an interval \( \omega \) and an escape time~\( N \).  A first obvious quantity of interest is the accumulated derivative along the critical orbit. We will compute  bounds for this and thus introduce the following notation: 
\begin{equation}
\label{eq:fnprime}
(f^{n})'(\omega) :=
\left[
\inf_{a\in \omega}(f^{n}_{a})'(c_{0}(a)), \ 
\sup_{a\in \omega}(f^{n}_{a})'(c_{0}(a))
\right].
\end{equation}
Also of great interest is the way in which the iterate \( c_{n}(a) \) of the critical point depends on the parameter. To study this dependeance, by some slight abuse of notation, let 
\(  c_{n}\colon\omega \to \omega_{n}  
\)
denote the map
\(
a \mapsto c_{n}(a). 
\) 
The map~\(  c_{n}  \)  is smooth with respect to \( a \) because the family \(  f_{a}  \) depends smoothly on the parameter, and so we  let \(  c_{n}'(a)  \) denote the derivative of \(  c_{n}  \) \emph{with respect to the parameter} (which is crucial in the parameter exclusion argument). Then we let 
\begin{equation}\label{eq:cnprime}
c'_{n}(\omega):=\left[\inf_{a\in \omega}c_{n}'(a), \sup_{a\in\omega}c_{n}'(a)\right].
\end{equation}
Also of interest, for less obvious and more technical reasons, in the parameter exclusion arguments, is the ratio between the derivatives with respect to the parameter and with respect to the phase space variable. We will therefore also compute the following quantities:
\begin{equation}
\label{eq:cfnprime}
\frac{c_{n}'}{(f^{n})'}(\omega) :=
\left[
\inf_{a\in \omega}\frac{c_{n}'(a)}{(f^{n}_{a})'(c_{0}(a))},\ 
\sup_{a\in \omega}\frac{c_{n}'(a)}{(f^{n}_{a})'(c_{0}(a))}
\right].
\end{equation}
Notice that bounds for \eqref{eq:cfnprime} can be easily derived from bounds for \eqref{eq:fnprime} and \eqref{eq:cnprime} but these may be quite far from optimal as there is no reason a priori for the lower and upper bounds in \eqref{eq:fnprime} and \eqref{eq:cnprime} to be attained for the same parameters. We will therefore compute bounds for \eqref{eq:cfnprime} directly. 

% ====================================================

\section{The Results}
\label{sec:results}

We now present and discuss the data obtained by our computations. In subsections \ref{sec:stocreg} and~\ref{sec:basic} we give a short overview of the procedure for subdividing the parameter space \( \Omega \) into a potentially large number of smaller subintervals. Then in the remaining subsections we give the statistics of several measurements which we carry out for these intervals. The raw data generated by our software is available in~\cite{data20}.

% ====================================================

\subsection{Stochastic and regular intervals} 
\label{sec:stocreg}
One of the results of our computations consists of a finite partition \( \mathcal P \) of \( \Omega \) made up of almost  \emph{4~million} explicit subintervals of  \( \Omega \). We will write \( \mathcal P \) as the union of two disjoint subsets 
\[
 \mathcal P^{+}=\{\text{``stochastic'' intervals}\} 
 \quad \text{ and } \quad 
 \mathcal P^{-} =\{\text{``regular'' intervals}\}
 \]
according to some rigorous and computationally verifiable properties of each interval. The construction of the partition \( \mathcal P \) depends on certain parameters of which the most important is the constant \( \delta>0 \) which defines the critical neighbourhood \( \Delta \), see \eqref{eq:Delta}. Recalling the notion of escape time in \eqref{eq:escape}, given any  \( N_{0}\geq 1 \), we will construct the collection of intervals \( \mathcal P^{+} \) so that 
 \[
 \emph{for each interval \( \omega\in\mathcal P^{+} \) there exists an escape time \( N\geq N_{0} \) for \( \omega \)}. 
 \]
The collection \( \mathcal P^{-} \) then consists simply of  intervals  for which the existence of such an escape time cannot be verified or, for whatever reason, is not  verified in our computations. 

The terminology ``regular interval'' and ``stochastic interval'' is only heuristic but suggestive of the fact that, while it is beyond the scope of this paper to prove this, it is reasonable to expect that most parameters in regular intervals are regular and most parameters in stochastic intervals are stochastic, as defined in \eqref{eq:regstoc}. For regular intervals this expectation is based on the data we compute, see  discussion at the end of  Section \ref{sec:exclusion}. For stochastic intervals this expectation is based on the arguments \cite{LuzTak06}, see discussion in Section \ref{sec:compstart}, and its verification is work in progress.

The purpose of this section is to describe the structure and properties of  \( \mathcal P^{+}\) and \(  \mathcal P^{-} \) for a particular choice of \( \delta \) and \( N_{0} \), namely 
\[
\delta=10^{-3} \quad \text{ and } \quad N_{0}=25.
\]
This particular choice of values is just for definiteness and does not have a particular meaning. 
Our main goal  is to show the kind of information that can be obtained by our computations. In particular, we will give rigorous estimates for the total measure of intervals in \( \mathcal P^{+}\) and \( \mathcal P^{-} \), as well as information about the distribution of the sizes of intervals, the computed values of \( N \), the sizes of \( \omega_{N} \), and other interesting information. 
The computations could just as well be carried out for any other values of \( \delta \) and \( N_{0} \), though they are clearly  more intensive and ``expensive'' for smaller values of \( \delta \) and larger values of \( N_{0} \). 

% ====================================================

\subsection{Basic strategy}
\label{sec:basic}

An important part of our approach  is that \( \mathcal P^{+}\) and \(  \mathcal P^{-} \) are \emph{dynamically defined}. We do not just try to verify the escape time condition in  some a priori given subdivision of \( \Omega \), but rather use dynamical information to subdivide  the parameter space \( \Omega \) in an efficient way. This makes a significant difference in terms of maximising the measure of intervals in \( \mathcal P^{+} \) and obtaining much more meaningful results. 
We describe this construction in detail in Section \ref{sec:P}. There are several non-trivial technical aspects to be addressed, especially in order to guarantee that all our estimates are rigorous, but the general strategy is actually very simple and we sketch it here. 

We start with the entire parameter space \( \omega=\Omega \) and consider the iterates \( \omega_{i} \) until they hit~\( \Delta \) at some time \( n\geq 1 \) (that is, \( \omega_n \cap \Delta \neq \emptyset \)). Then  we chop \( \omega \) into (at most 3) closed subintervals \( \omega=\omega^{\ell}\cup\omega^{\Delta}\cup\omega^{r} \) (the \emph{left}, \emph{middle}, and \emph{right} parts) with disjoint interiors in such a way that \( \omega^{\ell}_{n}\cap \Delta=\emptyset \) and \( \omega^{r}_{n}\cap \Delta=\emptyset \), and thus \( \omega_n \cap \Delta \subset \omega^\Delta_n \). We let \( \omega^{\Delta}\in \mathcal P^{-} \) and no longer consider any of its further iterations. If \( \omega^{\ell} \) is too small (according to some criteria specified precisely in Section~\ref{sec:assigning}), we also let it belong to  \( \mathcal P^{-} \) and stop iterating; the same with \( \omega^{r} \). Otherwise, we continue iterating \( \omega^{l} \) and \( \omega^{r} \) until they hit \( \Delta  \), and then we repeat the procedure. Every time an interval hits \( \Delta \), we verify whether \( n \geq N_{0} \) and the escape time condition holds. If this happens then we let the interval belong to \(  \mathcal P^{+} \) and stop iterating this interval. Moreover, if it is detected at any time during the computation of the iterates \( \omega_i \) that certain other conditions are met which suggest that none of further iterates of \( \omega \) is likely to lead to an escape time, we let the interval belong to~\( \mathcal P^{-} \) and stop iterating.

It is clear from the description of the construction that the collections \( \mathcal P^{+} \) and \( \mathcal P^{-} \) do not depend canonically on the choices of \( \delta \) and \( N_{0} \). Moreover, they also depend on some other choices; for example, on the level of binary precision \( p \) chosen for the computations, which we set as \( p := 250 \), on the minimum size \( w \)  (relative to \( \Omega \)) of an interval to consider it worth iterating, which we set as \( w := 10^{-10} \), and some other values relevant to the construction, as explained in detail in Section  \ref{sec:P}.  For the escape condition,we use the bound \( |\omega_N| \geq 0.0317 > \sqrt{\delta} \).  In a future paper we plan to analyse systematically the effect of changing these variables of the construction, but preliminary experiments indicate that while different choices may of course lead to quite different intervals being constructed, the overall statistics are remarkably stable and do not depend in a sensitive way on these choices, provided that we do not impose too severe restrictions on the computations, such as taking the precision \( p \) too low or the relative size \( w \) too large.

The computations were completed using the software described in Section~\ref{sec:software} and available at~\cite{software}. They were completed within 35 minutes  on a personal laptop computer with the Intel\textsuperscript{\textregistered} Core{\texttrademark} i5-8265U processor. The results of the computations are available in~\cite{data20}.

% ====================================================

\subsection{Measure of regular and stochastic intervals}
\label{sec:exclusion}

The partition \( \mathcal P \) obtained by our computations is made up of the disjoint union of the families of \( \mathcal P^{+} \) and \( \mathcal P^{-} \) made up respectively of stochastic intervals, which satisfy the escape time condition~\eqref{eq:escape}, and regular intervals, for which this condition was not  verified. The fundamental quantities of interest are therefore the number and total measure of the intervals in each family.  The first and most striking observation is that 
\begin{quote}
\emph{almost 90\% of parameters belong to stochastic intervals.} 
\end{quote}
This means that 90\% of parameters belong to intervals which have an escape at some relatively large time \( N\geq 25 \) (and are therefore good candidates for the parameter exclusion arguments).  More precisely, letting \( \#\mathcal P \) denote the cardinality of the partition \( \mathcal P \) and, by some slight abuse of notation, letting \( |\mathcal P| \) denote the total measure of intervals in \( \mathcal P \), we have the following results. The partition \( \mathcal P \) is formed by almost 4 million intervals or, more precisely, 
\[
\#\mathcal P  =  3,\!969,\!763 \quad \text{ and } \quad |\mathcal P|=0.6 = |\Omega|.
\]
Of these, about 36\% in number and 90\% in measure are stochastic, more precisely 
\[
\#\mathcal P^{+} = 1,\!436,\!063 \geq 0.36 \#\mathcal P
 \quad \text{ and } \quad
 |\mathcal P^{+}|\geq 0.539934844013 \geq 0.89989 |\Omega|, 
\]
and therefore 
\[
\#\mathcal P^{-} = 2,\!533,\!700\ \leq 0.64 \#\mathcal P
 \quad \text{ and } \quad
 |\mathcal P^{-}| \leq 0.060065155986 \leq 0.10011 |\Omega|. 
\]

In the following subsections we analyse in detail several properties of the family \( \mathcal P^{+} \) of stochastic intervals, which are our main objects of interest. It is worth, dwelling a little bit here on the collection \( \mathcal P^{-} \) of regular intervals, which also exhibit some very interesting features. First of all, as many intervals in \( \mathcal P^{-} \) are adjacent to each other (the same is also true in \( \mathcal P^{+} \)), it can be useful to merge adjacent intervals and consider ``connected components'' of \( \mathcal P^{-} \) which are a bit less dependent on the specifics of the construction. In terms of these connected components, it is interesting to observe that the total measure of \( \mathcal P^{-} \)is disproportionately concentrated on larger intervals. The 100 largest  components (actually made up of \( 1{,}124{,}307 \)  intervals of \( \mathcal P^{-} \)) take up a total measure of about  \( 0.05726 \), which is 95\% of the total measure of \( \mathcal P^{-} \), and the 3 largest components (made up of \( 11{,}830 \), \( 7{,}955 \) and \( 8{,}313 \) intervals respectively) alone take up more than 30\% of the total measure of \( \mathcal P^{-} \). These largest 3 are contained in the following intervals: 
\begin{eqnarray*}
&&I_{1} = [1.75208241722, 1.77992046728], \quad |I_{1}|= 0.0278381, \\
&&I_{2} = [1.47590994781, 1.48293277717], \quad  |I_{2}|= 0.00702283, \\
&&I_{3} = [1.62533272418, 1.63110961362],  \quad |I_{3}| = 0.00577689. 
\end{eqnarray*}   

\begin{figure}[htb]
\begin{center}
\includegraphics[width=\textwidth]{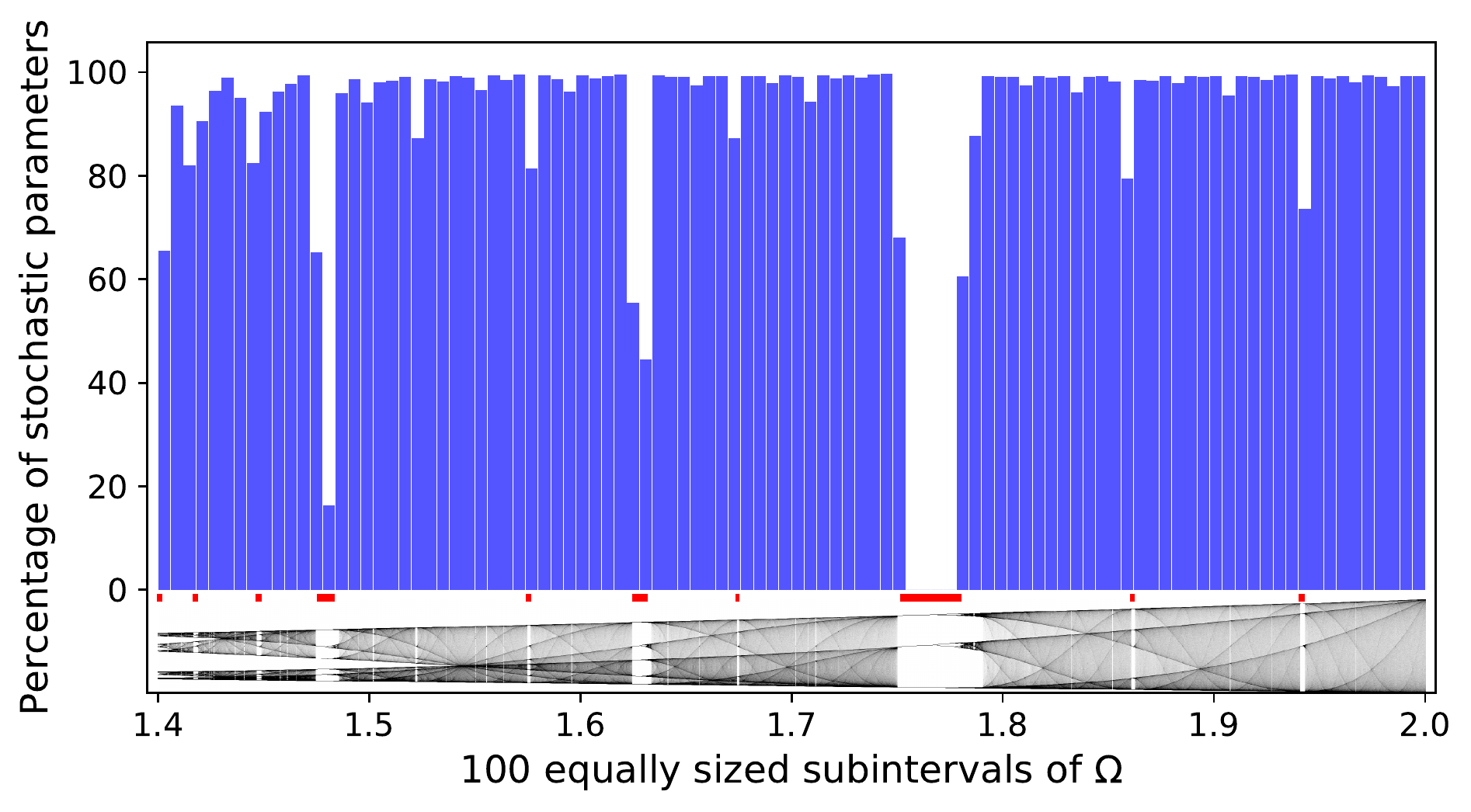}
\caption{Distribution of the measure of stochastic parameters in \( \Omega \).
Blue bars show the percentage of stochastic parameters in each of the 100 subintervals of \( \Omega = [1.4,2] \).
Horizontal red lines show the location of 10 largest connected components of \( \mathcal{P}^- \).
Bifurcation diagram for the quadratic map is shown along the horizontal axis. More detailed discussion of this picture at the end of Section \ref{sec:exclusion}.}
\label{fig:stochasticintervals}
\end{center}
\end{figure}

In Figure \ref{fig:stochasticintervals}, we have represented the ten largest connected components of \( \mathcal P^{-} \) by red horizontal bars to highlight how  they match up remarkably well, albeit unsurprisingly, with the well known \emph{periodic windows} which appear in the standard \emph{bifurcation diagram}. It would clearly be interesting to prove that most parameters in regular intervals are indeed regular, perhaps adapting the techniques of \cite{TucWil09}.
For clarity and completeness, we remark that Figure \ref{fig:stochasticintervals} was created by dividing  \( \Omega = [1.4,2] \) into \( 100 \) intervals \( \omega_1, \ldots, \omega_{100} \) of the same width \( 0.006 \). For each of these intervals \( \omega_i \), the corresponding blue bar shows the percentage of stochastic parameters in \( \omega_i \), that is, the measure of \( \omega_i \cap P^+ \), where \( P^+ \subset \Omega \) is the union of all the intervals in \( \mathcal{P}^+ \). The height of a blue bar below 100\% indicates that \( \omega_i \) intersects some intervals in \( \mathcal{P}^- \). Note that the height of the bars above a large periodic window in the bifurcation diagram (shown along the horizontal axis) is zero if the corresponding \( \omega_i \) is entirely covered by intervals in \( \mathcal{P}^- \). The alignment of bars considerably lower than 100\% with the periodic windows clearly shows how the periodic windows contribute to \( \mathcal{P}^-\). For example, from the graph (or actually from raw data that was used to plot the graph) one can read that about 99.34\% of the interval \( [1.988,1.994] \) is covered by \( \mathcal{P}^+ \), while only some 73.6\% of the interval \( [1.94,1.946] \) is covered by \( \mathcal{P}^+ \).

% ====================================================

\subsection{Distribution of sizes of stochastic intervals}

\begin{figure}[htbp]
\begin{center}
\includegraphics[width=0.7\textwidth]{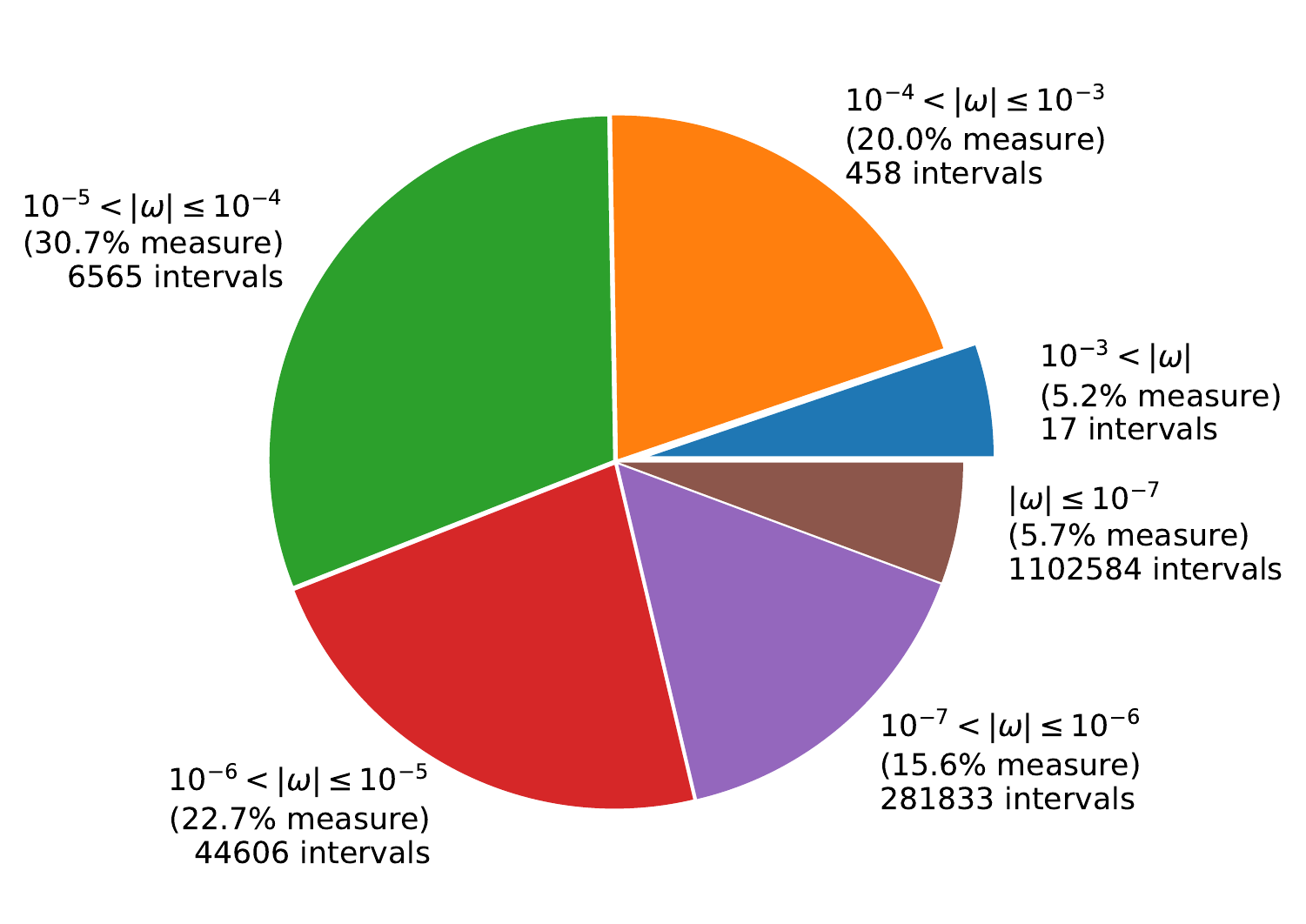}
\caption{Distribution and number of stochastic intervals of different sizes}
\label{fig:stochdist}  
\end{center}
\end{figure}

We now focus on the family \( \mathcal P^{+} \) of stochastic intervals, which occupy almost 90\% of the parameter space \( \Omega \) and are our main objects of interest. 
Figure \ref{fig:stochdist}  shows  the distribution of sizes of stochastic intervals, which turns out to span several orders of magnitude of different scales. Most of the measure is taken up by ``medium'' to ``small'' intervals, whereas  ``large'' intervals (\( |\omega|\geq 10^{-3} \)), and ``very small'' intervals (\( |\omega|\leq 10^{-7} \)) each take up about 5\% of the total measure.  Notice  that, perhaps also unsurprisingly, the number of very  small intervals is more than the number of intervals of all other sizes put together. This seems to suggest that, similarly to the regular intervals, while the number of small intervals grows quite fast, it does not grow fast enough to have a significant effect on the measure. 

% ====================================================

\subsection{Distribution of escape times \( N \)}

\begin{figure}[htbp]
\begin{center}
\includegraphics[width=0.75\textwidth]{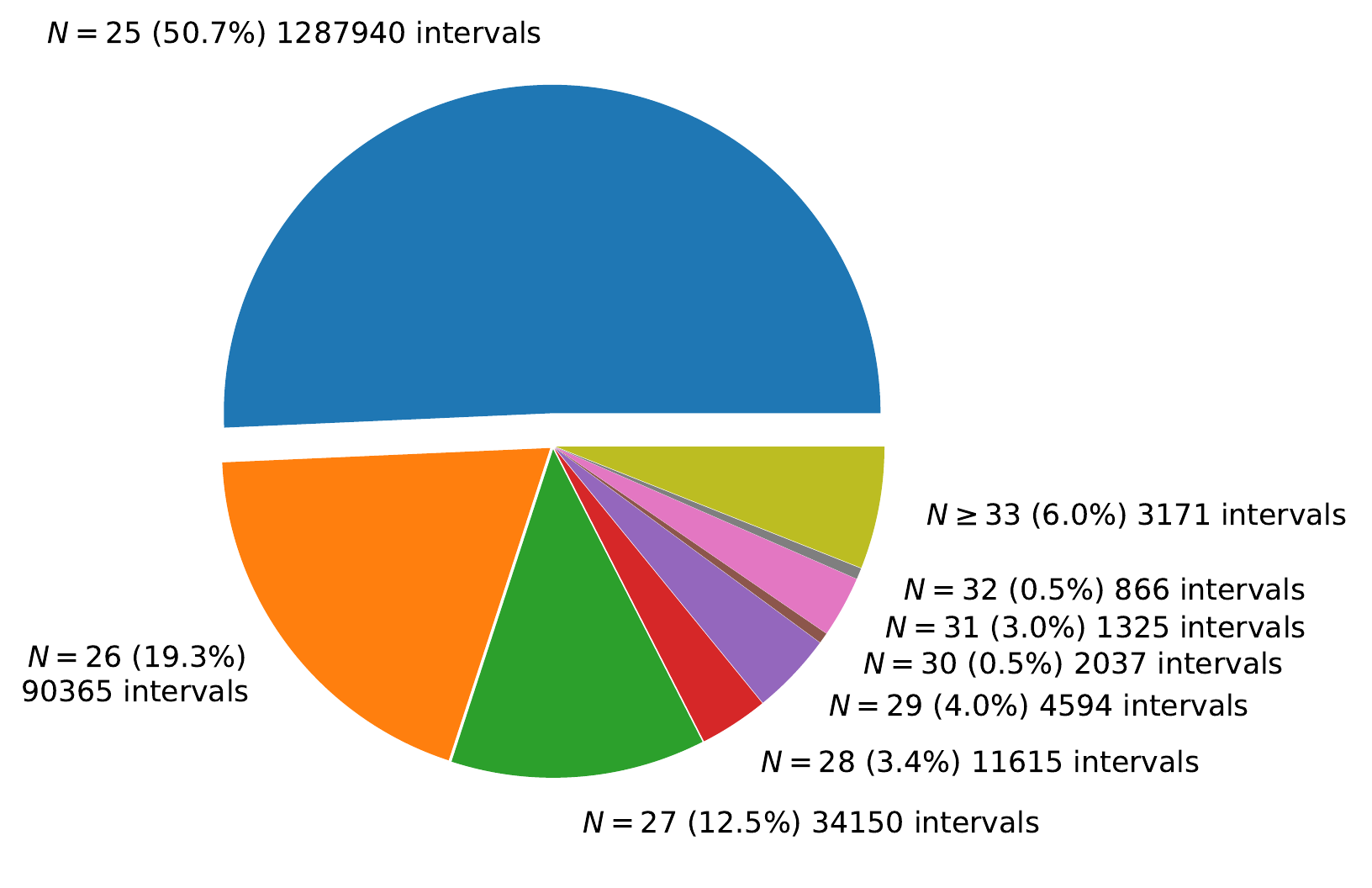}
\caption{Distribution of escape times of stochastic intervals}
\label{fig:distescape}  
\end{center}
\end{figure}
 
Figure \ref{fig:distescape}  shows the distribution of escape times of stochastic intervals. 
Remarkably, more than 90\% of the intervals, occupying more than 50\% of the measure,  escape at the very first opportunity, with escape time \( N=N_{0}=25 \). Most other intervals have escape times just slightly larger than 25, with more than 99.7\% of intervals, occupying 94\% of the measure, having escape times \( 25\leq N \leq 32 \). We, emphasize, however that there is a long ``tail,'' and intervals exist with much higher escape times, up to a maximum of escape time \( N=199 \) for 73 distinct intervals in \( \mathcal P^{+} \) taking up 0.00173\% of the total measure of stochastic intervals. 

% ====================================================

\subsection{Distribution of sizes of intervals \( \omega_{N} \) at escape times}

\begin{figure}[htbp]
\begin{center}
\includegraphics[width=0.75\textwidth]{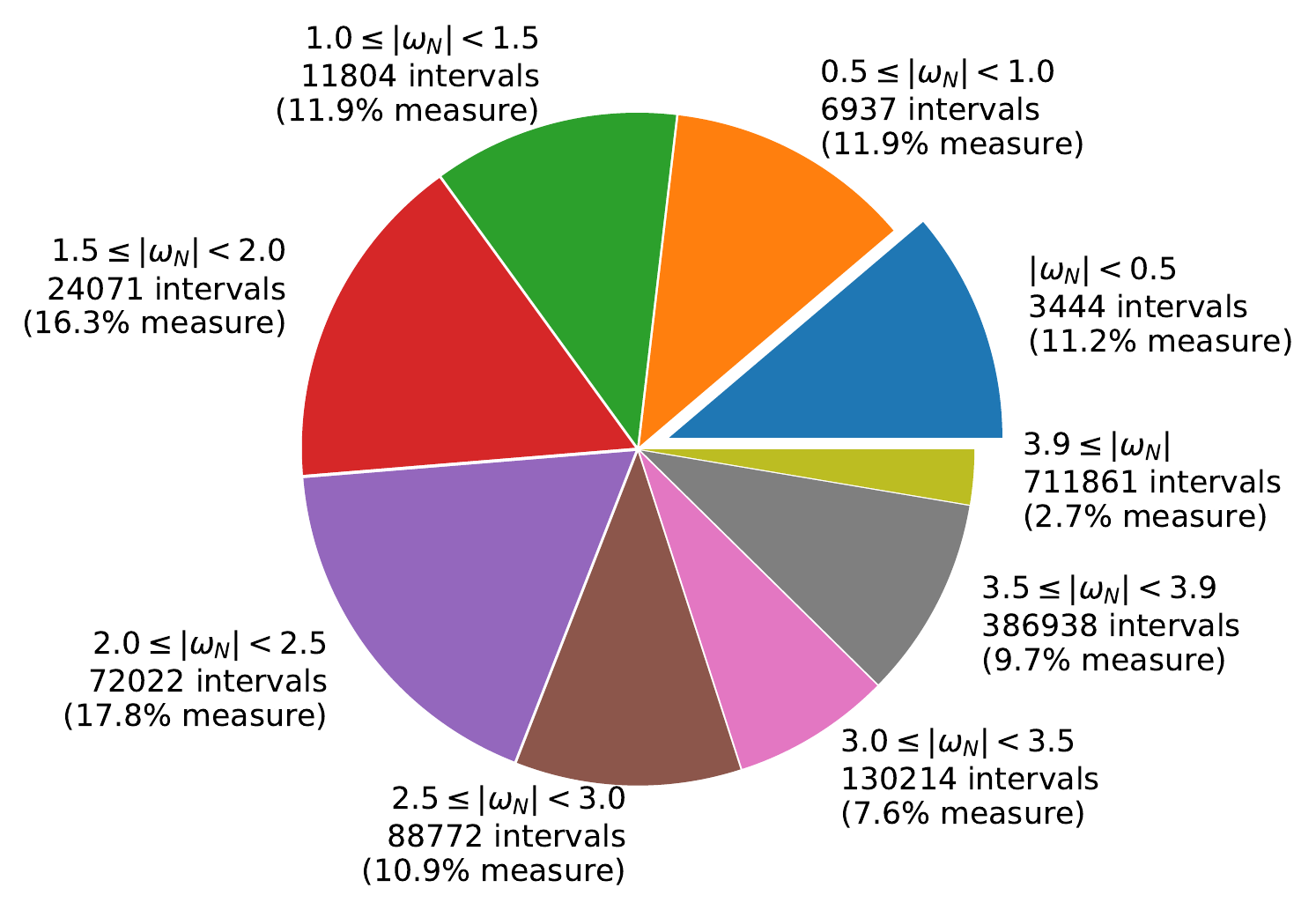}
\caption{Distribution of sizes of intervals at escape times}
\label{fig:excludedmeasure}  
\end{center}
\end{figure}

Figure \ref{fig:excludedmeasure} shows another distribution, namely the sizes of the images \( \omega_{N} \) of stochastic intervals  at their escape times. The results  are, in our opinion, quite \emph{unexpected and interesting} even though, given the unpredictable way intervals are regularly chopped as part of the construction of the partition \( \mathcal P \), there seems to be no elementary heuristic argument for predicting the size of \( \omega_{N} \). 
Recall that by definition of escape time we always have  a lower bound of \(  0.0317  \gtrsim \sqrt \delta\), which is relatively small in relation to the interval \( I_{a} \)  of definition of the map (which is 4 for the ``top'' parameter \( a=2 \) and slightly less for other parameters). It seems  therefore quite remarkable that more than 99.7\% of intervals, occupying almost 90\% of the measure, have relatively ``macroscopic'' size, with \( |\omega_{N}| \geq 0.5 \). Even more, it turns out that intervals occupying some 20\% of the measure have ``very large'' images,  i.e.  \( |\omega_{N}|\geq 3 \). In Section \ref{sec:study} we analyse in detail the ``personal history'' of one, more or less randomly chosen, interval \( \omega\in\mathcal P^{+} \) with escape time \( N=26 \) and such that \( |\omega_{N}|\geq 3.5 \),  in order to help understand the mechanism by which this situation can occur. 

Figure~\ref{fig:excludedmeasure} reveals one more interesting piece of information.
The first few pieces in the pie chart show the measure of intervals \(\omega\) that yield smaller \( \omega_{N} \).
Although their measure is considerable, the actual number of intervals that yield this measure is not very big. For example, the first 4 pieces that yield over 50\% of the measure consist of only 44,256 individual intervals.
On the other hand, the last two pieces of the pie chart, corresponding to the largest \( |\omega_{N}| \), comprise as little as 12.4\% of the measure, yet they consist of almost 1.1 million intervals.
This shows that there are many large intervals that yield small \( \omega_N \)  and many tiny intervals that yield huge \(\omega_N\); one could call it \emph{negative correlation} between \( |\omega| \) and \( |\omega_N| \).

% ====================================================

\subsection{Accumulation of derivatives}
\label{sec:accum}

\begin{figure}[htbp]
\begin{center}
\includegraphics[width=0.75\textwidth]{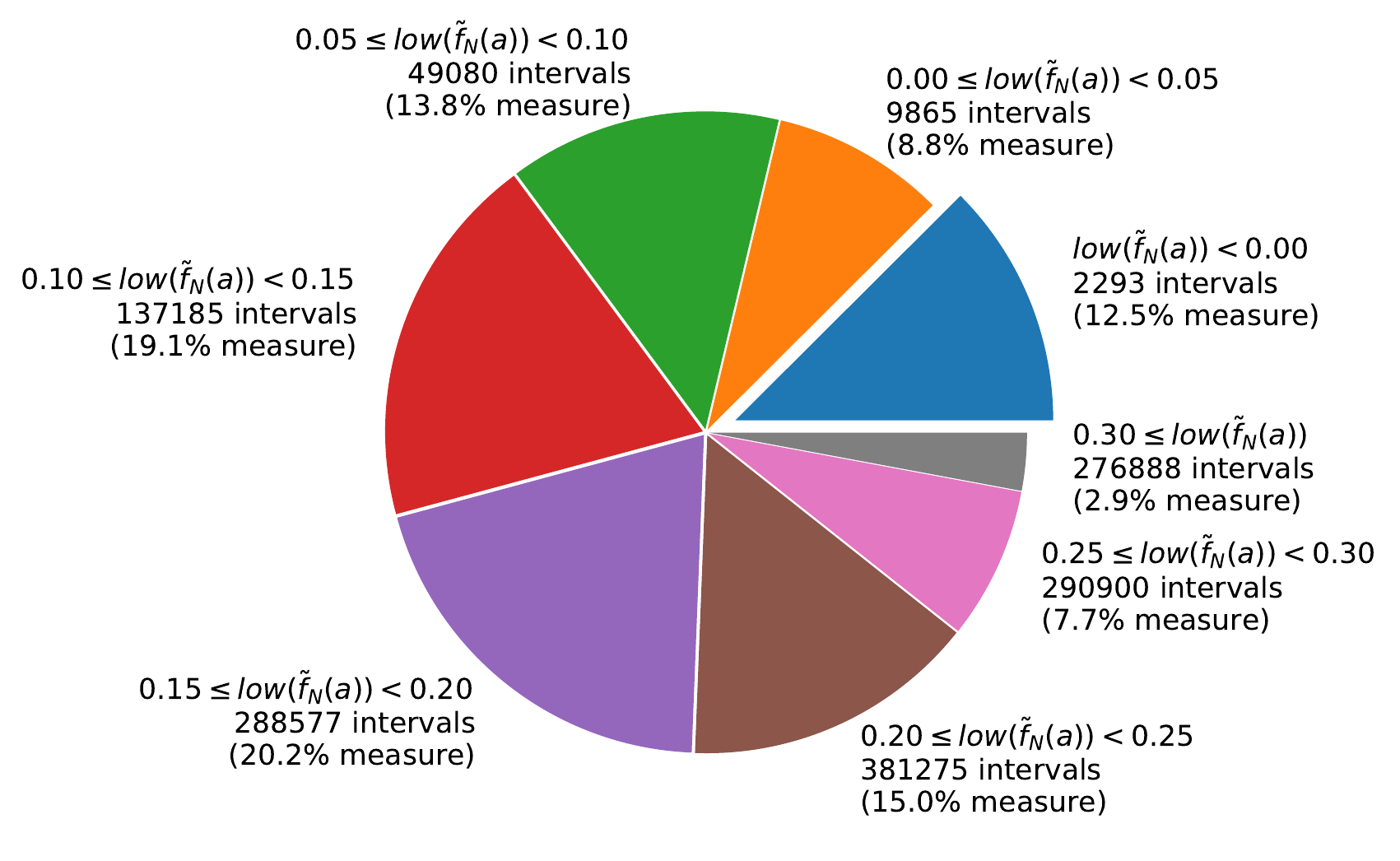}
\caption{Distribution of the lower bounds on \( \tilde f_{N}(\omega) \) computed for the stochastic intervals.
The highest encountered value was \( \approx 0.732 \).}
\label{fig:rateF}
\end{center}
\end{figure}

\begin{figure}[htbp]
\begin{center}
\includegraphics[width=0.75\textwidth]{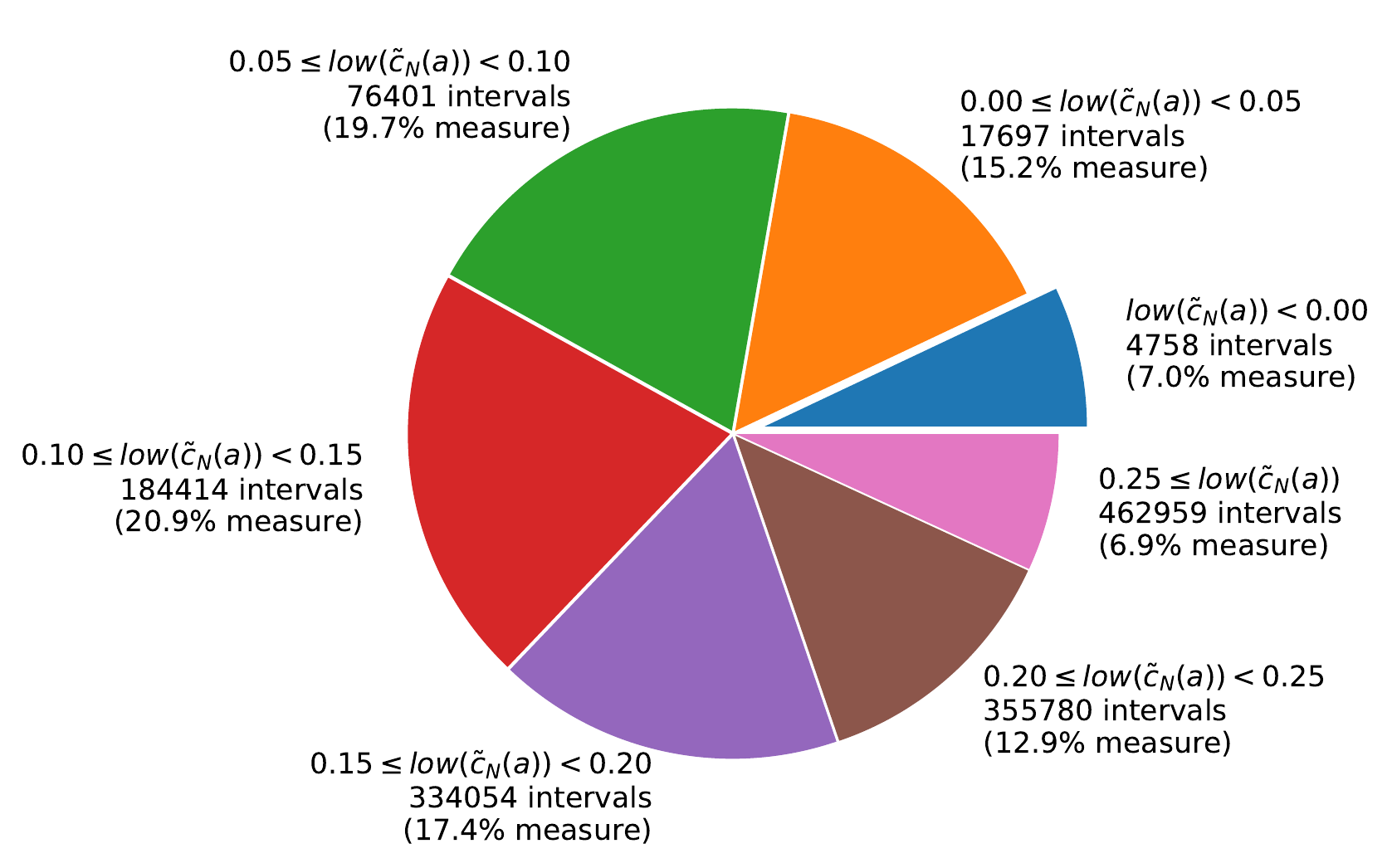}
\caption{Distribution of the lower bounds on \( \tilde c_{N}(\omega) \) computed for the stochastic intervals.
The highest encountered value was \( \approx 0.716 \).}
\label{fig:rateC}
\end{center}
\end{figure}

Figures \ref{fig:rateF} and \ref{fig:rateC} show some results related to the computations of the space and parameter derivatives, as in  \eqref{eq:fnprime} and \eqref{eq:cnprime},  These can be of significant interest in a variety of contexts, especially when they exhibit \emph{exponential growth}, which is a non-trivial feature, given that some iterates can be very close to the critical point where the derivative vanishes. In view of this fact, and of the large variation in the escape times for stochastic intervals, it seems best to present the data in the form of average \emph{exponential rate} of growth along the orbits. Thus, for a stochastic interval \( \omega\in \mathcal P^{+} \) with escape time \( N \), and  a parameter \( a\in \omega \), we define 
\[
\tilde f_{N}(a):= \frac 1N \log |(f^{n}_{a})'(c_{0}(a))|
\quad\text{ and } \quad 
\tilde c_{N}(a):= \frac 1N \log |c'_{N}(a)|
\]
and then, analogously to \eqref{eq:fnprime} and \eqref{eq:cnprime}, we define 
\[
\tilde f_{N}(\omega) :=
\left[
\inf_{a\in \omega} \tilde f_{N}(a), \ 
\sup_{a\in \omega} \tilde f_{N}(a) 
\right] 
\quad\text{ and } \quad 
\tilde c_{N}(\omega) :=
\left[
\inf_{a\in \omega} \tilde c_{N}(a), \ 
\sup_{a\in \omega} \tilde c_{N}(a) 
\right].
\]

Figure \ref{fig:rateF} shows the distribution of lower bounds computed for \(  \tilde f_{N}(\omega) \). We note that for 12.5\% of intervals in measure we do not have a positive lower bound, but this does not necessarily mean that there is no exponential growth. Indeed, all these intervals have a positive upper bound (not represented here) and it seems most likely that the lack of a positive lower bound is due to overestimates caused by using interval arithmetic in evaluating these quantities. We also note that there is a remarkably even distribution of lower bounds, with about 10-20\% in measure of parameter intervals in each band, except for the highest rates of growth above \( 0.3 \) which is exhibited only by 2.9\% of parameters. We mention, however, that higher rates are exhibited by smaller fractions of parameters, all the way up to 0.732. 

\begin{figure}[tb]
\begin{center}
\includegraphics[width=0.8\textwidth]{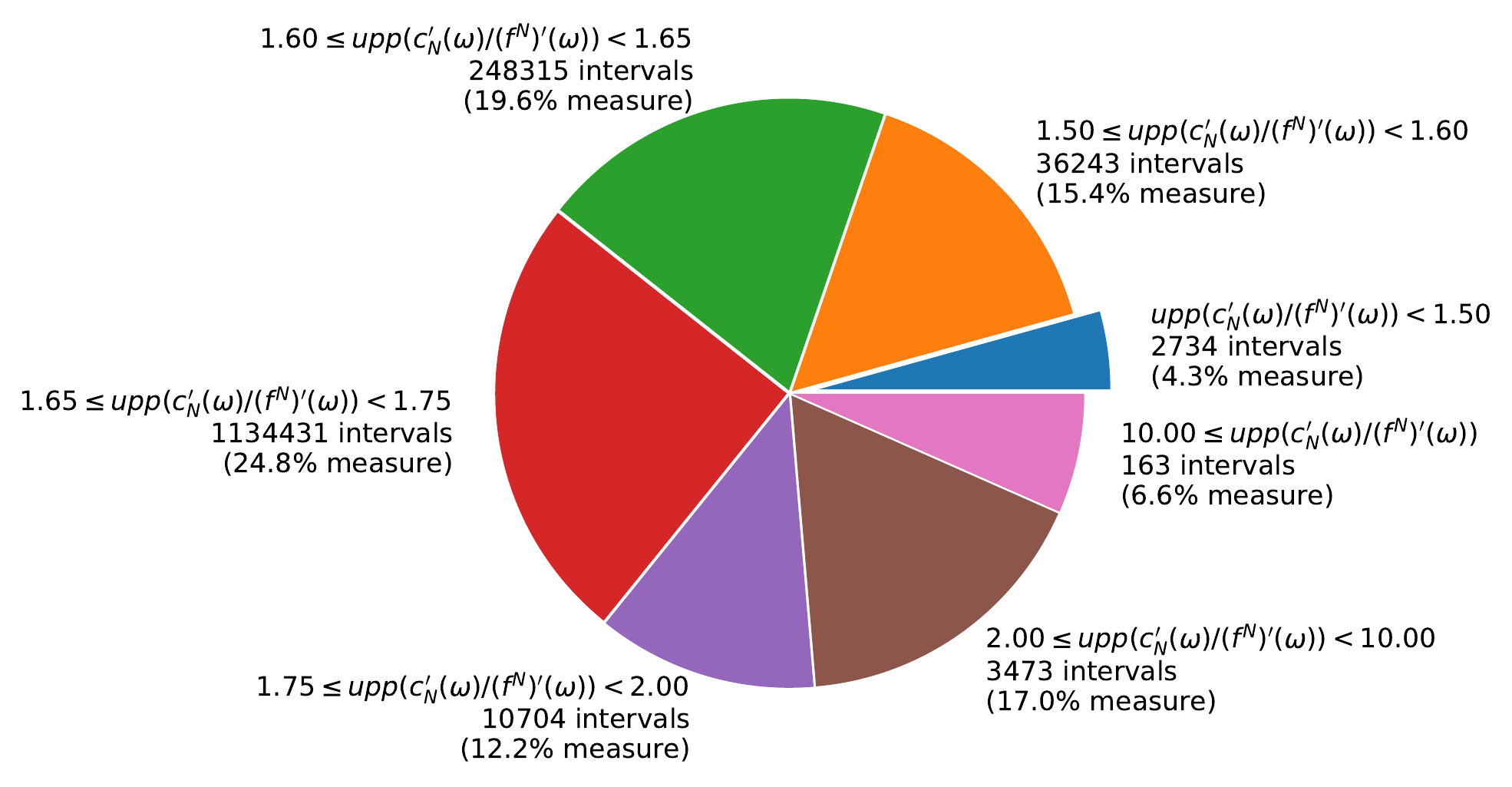}
\caption{Distribution of the upper bounds on the quotient of derivatives \( c_{N}'(a) / (f^{N}_{a})'(c_{0}(a)) \) computed for the stochastic intervals.
The lowest encountered value was almost \( 1.2 \), the highest was close to \( 10^{13}  \).}
\label{fig:D}
\end{center}
\end{figure}

\begin{figure}[htbp]
\begin{center}
\includegraphics[width=0.75\textwidth]{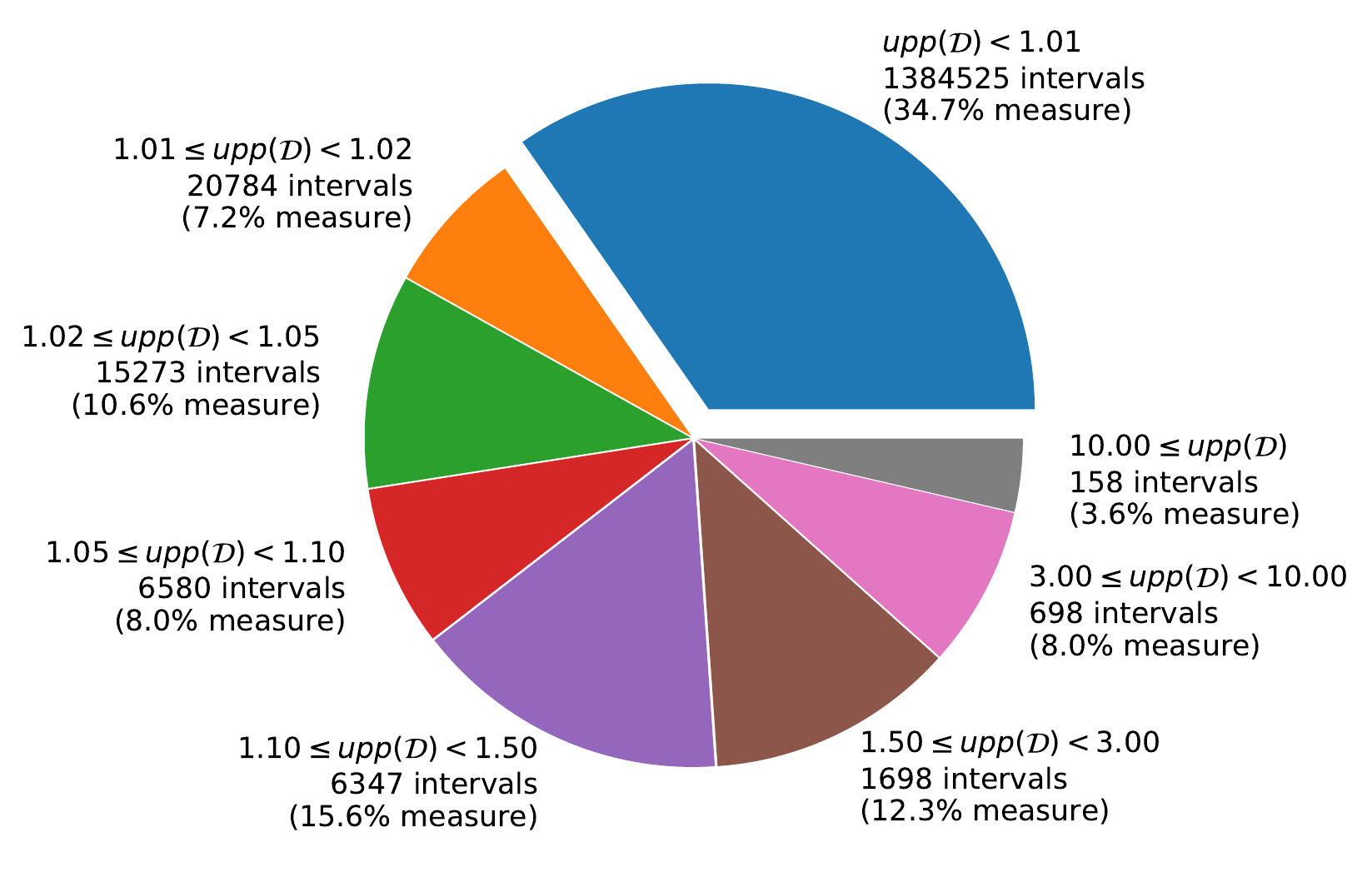}
\caption{Distribution of the upper bounds on the quotient \( \mathcal{D} \) computed for the stochastic intervals.
The lowest encountered value was slightly above \( 1 \), the highest was almost \( 400 \).}
\label{fig:quotD}
\end{center}
\end{figure}

Figure \ref{fig:rateC} shows the corresponding statistics for \( \tilde c_{N}(\omega) \) which turn out to be remarkably similar to those for \(  \tilde f_{N}(\omega) \). We note however that the close relationships between these values is ``real'', not just statistical, as demonstrated in Figures \ref{fig:D} and \ref{fig:quotD} which refer to the measurements of the ratio \eqref{eq:cfnprime} between these two quantities. 
 Figure \ref{fig:D} shows the statistics of the upper bounds for this ratio, and should be interpreted  in conjunction with Figure \ref{fig:quotD}  which gives upper bounds for the \emph{distortion} 
 \[
   \mathcal D:= \frac{\sup_{a\in\omega}\{|c'_{N}(a)/(f^{N}_{a})'(c_{0}(a))|\} }
   {\inf_{a\in\omega}\{|c'_{N}(a)/(f^{N}_{a})'(c_{0}(a))|\}}.
\]
It seems highly remarkable that this distortion is very close to 1 in most of the intervals, both in cardinality and in measure, and \( <1.5  \) for more than 75\% of intervals, both in cardinality and in measure. This means that for most parameters the upper and lower bounds for \( |c'_{N}(\omega)/(f^{N}_{a})'(c_{0}(\omega))| \) are comparable and thus Figure \ref{fig:D} gives a good representation of its actual values. 
It seems therefore also highly remarkable that this ratio is  \( < 2 \) for over 75\% in the measure of parameter intervals. 

% ====================================================

\section{The Computations}
\label{sec:strategies}

In order to cater for readers with different levels of familiarity with computational methods, in Sections \ref{sec:strategies}-\ref{sec:algorithms} we give increasingly detailed and technical description of the computational procedures and algorithms used to obtain the results given in Section \ref{sec:results}. We begin, in this section, by explaining our general strategy for the computations, in a way that is easily accessible to anyone with some familiarity with one-dimensional dynamics, emphasising nevertheless some crucial but subtle aspects related to the need to obtain rigorous explicit bounds. In Section \ref{sec:P} we give a detailed but non-technical explanation of the procedure for constructing the families of intervals \( \mathcal P^{-} \) and \( \mathcal P^{+} \) using the results of the calculations described in this section.  In Section \ref{sec:procedure} we explain how the calculations can be  formalised in order to work with computer representable numbers and to yield rigorous bounds for all the quantities we compute. Finally, in Section \ref{sec:algorithms} we describe precisely the  algorithms used to implement each step of the procedure. 

The computations can be divided roughly into three categories, which  we describe  in the following three subsections. 
 
% ====================================================

\subsection{Iterating}
\label{sec:iterate}
The core challenge we address in this paper is the development of  effective techniques for the computation of intervals   \( \omega_{n}:=\{c_{n}(a): a\in \omega\} \), for some given parameter interval \( \omega\subseteq \Omega \). 

The first step in this direction is clearly the development of effective techniques for the computation of the point \( c_{n}(a):=f^{n}_{a}(c_{0}(a))=f^{n}_{a}(a) \) for a fixed parameter  \(  a\in \omega  \) (recall that \( c=0 \) and so \( c_{0}(a)=f_{a}(c)=a \)). This is already non-trivial since the value of \( a \) may not be computer representable and therefore require an approximation strategy before we even begin iterating. Even if \( a \) is representable, its first image \( c_{1}(a):=f_{a}(a)=a-a^{2} \) is very possibly not representable, and similarly for higher iterates. Fortunately, tried and tested methods, known as  \emph{interval arithmetic} \cite{Moore1966, WT2011}, exist and can be very effective for these kinds of computations. They consist essentially of \emph{enclosing} the point to be iterated in a small interval whose endpoints are representable numbers, and then applying the map to this interval to obtain a rigorous \emph{enclosure}, and therefore an approximation, of the image of the given point. The method of course gives increasingly large enclosures, and therefore increasingly poor approximations, for higher iterates \( c_{n}(a) \) but these can still be obtained to any desired precision for a fixed \( n \) by increasing the computer precision and therefore the cardinality, and ``density'', of the set of representable numbers. For example, in the calculations in Section \ref{sec:results} we work with about 80  decimal places. 

In principle we  could blindly apply the interval arithmetic techniques also to the computation of the intervals \( \omega_{n} \). Indeed, supposing for example that the parameter interval \( \omega=[\mathbf{a},\mathbf{b}] \) was given by endpoints which are representable numbers (here and below we will conventionally use bold type to denote representable numbers) and that the same was true of the interval \( \omega_{i} := [\mathbf{a}_{i},\mathbf{b}_{i}] \) for some \( i\geq 0 \) (or that we had a representable enclosure of the interval \( \omega_{i} \), this does not make much of a difference for the discussion here). Then we could use interval arithmetic to compute a rigorous enclosure for all possible values of \( a-x^{2} \) for all possible \( a\in \omega \) and \( x\in \omega_{i} \), thus yielding a rigorous enclosure for \( \omega_{i+1} \). It is easy to see, however, that this will very likely produce \emph{huge} overestimates of \( \omega_{i+1} \), which would moreover compound at each iteration, and is therefore not at all a very \emph{effective} way to proceed.  The reason for the overestimation is due to the fact that this approach consists of iterating \emph{every} point in \( \omega_{i} \) by \( f_{a} \) for \emph{every} parameter \( a\in \omega \), rather than iterating each point in \( \omega_{i} \)  just by the corresponding parameter. The enclosure for \( \omega_{i+1} \) will therefore contain the points \( f_{\mathbf{a}}(\mathbf{b}_{i}) \) and \(  f_{\mathbf{b}}(\mathbf{a}_{i})  \) which may be much further apart than necessary if, for example,  \( \omega_{i} \) lies on the right of the critical point. 

At first sight, there is an obvious solution to this problem, which is to simply iterate the points corresponding to the endpoints \( \mathbf{a} \) and \( \mathbf{b} \) of the parameter interval \( \omega \), i.e., to compute the points \( c_{n}(\mathbf{a}) \) and \( c_{n}(\mathbf{b})  \). As mentioned above, the computation of these points can be easily achieved to arbitrary precision. The problem, however, is that it is  not necessarily the case that \( c_{n}(\mathbf{a}) \) and \( c_{n}(\mathbf{b})  \) are the endpoints of \( \omega_{n} \) even though \( \mathbf{a} \) and \( \mathbf{b}  \) are the endpoints of \( \omega \), since  the map \( c_{n}\colon\omega\to\omega_{n} \) may fail to be injective, and if it is not injective then it may ``fold'' and one of \( c_{n}(\mathbf{a}) \) or \( c_{n}(\mathbf{b})  \) may lie in the interior of \( \omega_{n} \). We can resolve this issue if, recalling \eqref{eq:cnprime}, we have 
\begin{equation}\label{eq:mon}
0\notin c_{n}'(\omega)
\end{equation}
which implies that  \(  c_{n}'(a) \neq 0  \) for every \(  a\in \omega  \) and therefore  that the map \(  c_{n}  \) is \emph{monotone} on \(  \omega  \).  This  implies that \( c_{n}(\mathbf{a}) \) and  \( c_{n}(\mathbf{b})  \) are indeed the endpoints of \(  \omega_{n}  \) and therefore provides both  \emph{inner} and \emph{outer} enclosures of \( \omega_{n} \) to arbitrary precision. 

We emphasize that \eqref{eq:mon} cannot always be verified and that its verification is implicitly one of the conditions required for a  parameter interval \( \omega \) to belong to \( \mathcal P^{+} \). As mentioned above, the collection \( \mathcal P^{-} \) is formed by those intervals for which the escape condition cannot be verified, for a variety of possible reasons, and failure to satisfy \eqref{eq:mon} is one of these reasons. 
We will describe below the precise way in which we check \eqref{eq:mon}, we just mention here that it will be done by a simple inductive procedure. For \( n=0 \), we have \( c_{0}(a)=a \), and therefore \( c_{0}'(a) = 1 \) for all \(a \in \omega\). For  \( n\geq 1 \), we use the formula
\begin{equation}\label{eq:derivC}
c_n'(a)=-2c_{n-1}(a) \cdot c_{n-1}'(a)+ 1.
\end{equation}
%so that
%\(
%c_n'=(-2)^{n}c_{n-1}\cdots c_{0}+(-2)^{n-1}c_{n-1}\cdots c_{1}+\cdots+(-2)^1c_{n-1}+1,
%\)
If we have rigorous enclosures for both \(  \omega_{n-1}  \) and \(  c'(\omega_{n-1})  \) then we can use \eqref{eq:derivC} and standard interval arithmetic computations to obtain a rigorous enclosure for \(  c'_{n}(\omega)  \) and check  \eqref{eq:mon}. 

% ====================================================

\subsection{Differentiating}
\label{sec:estimate}

As mentioned in the introduction, we are also interested in computing rigorous enclosures for the intervals \( (f^{n})'(\omega) \) and \( c_{n}'/(f^{n})'(\omega) \) defined in \eqref{eq:fnprime} and \eqref{eq:cfnprime}.
For both intervals we use an inductive procedure similar to that used for the calculation of \(  c'_{n}(\omega)  \) above. Specifically, by the chain rule we have 
\begin{equation}\label{eq:chain}
\begin{aligned}
(f^{n}_{a})'(c_0(a)) &= f'_{a}(f^{n-1}_{a}(c_0(a))) (f^{n-1}_{a})'(c_0(a)) 
\\ &= -2 f^{n-1}_{a}(c_0(a)) (f^{n-1}_{a})'(c_0(a)) 
\\ &= -2 c_{n-1}(a) (f^{n-1}_{a})'(c_0(a))
\end{aligned}
\end{equation}
and therefore, using interval arithmetic, rigorous enclosures for \(  \omega_{n-1}  \) and for \(  (f^{n-1})'(\omega)  \) immediately yield rigorous enclosures for \(  (f^{n})'(\omega)  \). 
Similarly, \eqref{eq:derivC} and \eqref{eq:chain} imply 
\begin{equation}\label{eq:quot}
\frac{c'_{n}(a)}{(f^{n}_a)'(c_0(a))}=  \frac{c'_{n-1}(a)}{(f^{n-1}_a)'(c_0(a))}+\frac{1}{(f^{n-1}_a)'(c_0(a))}
\end{equation}
and therefore, rigorous enclosures for \(  c'_{n-1}/(f^{n-1})'(\omega)  \) and for \(  (f^{n-1})'(\omega)  \) yield a rigorous enclosure for  \(  c'_{n}/(f^{n})'(\omega)  \). 
Notice that an enclosure for \(  c'_{n}/(f^{n})'(\omega)  \) could also be  computed directly from the enclosures of \(  c'_{n}  \) and~\(  (f^{n})'  \) by taking the worst case bounds, but the bounds we compute here,  using \eqref{eq:quot} inductively, are clearly much sharper.

% ====================================================

\subsection{Chopping}
\label{sec:avoid}
Finally, we discuss in a bit more detail, the chopping procedure described briefly at the end of Section \ref{sec:stocreg} leading the the construction of the partition \( \mathcal P \). As described there, the basic strategy is very simple and intuitive, chopping intervals which hit the critical neighbourhood \( \Delta \), say at some time \( n\geq 1 \), into subintervals which either land outside \( \Delta  \) and can be iterated further, or continue to intersect  \( \Delta \) and therefore belong to \( \mathcal P^{-} \). The computational problem  is simply stated and consists of finding the boundary between the parameters which fall into \( \Delta \) and those which do not. While we can have a very good approximation of the entire interval \( \omega_{n} \) following the procedure described in Section \ref{sec:iterate}, this is based on computation of the endpoints and does not help in identifying the parameters in the interior of \( \omega \) whose images fall into a particular position, such as close to the boundary of \( \Delta \). The map \( c_{n} \colon \omega\to \omega_{n}\) is not affine and therefore we cannot directly recover the parameters in \( \omega \) which map to the boundary points of \( \Delta \) under \( c_{n} \), even knowing with a good degree of accuracy the position of the boundary points of \( \omega_{n} \). 

Our approach is to use a relatively straightforward variant of the numerical algorithm known as the \emph{bisection method} (see e.g. \cite[\S 3.1]{KinChe1991}).  In order to explain this approach, let us fix \( \omega = [u,v] \) and \( n > 0 \).
Assume condition \eqref{eq:mon} holds true, and \( \omega_n \cap \Delta \neq \emptyset \).
Assume \( c_n (u) \notin \Delta \).
For simplicity of notation, assume \( c_n \) is increasing.
To make the idea clear, let us ignore rounding errors for the moment and assume the computations are exact.
We are going to construct inductively two sequences of numbers \( \{x_i\} \) and \( \{y_i\} \), with \( x_i < y_i \)
and \( |y_i - x_i| = 2^{-i} |v - u| \),
with the following property: \( c_n ([u,x_i]) \cap \Delta = \emptyset \) and \( c_n ([u,y_i]) \cap \Delta \neq \emptyset \).
In this way, by computing consecutive elements of the two sequences, we are going to get a gradually better approximation
of \( c_n^{-1} (-\delta) \).
Set \( x_0 := u \) and \( y_0 := v \), which satisfies the required properties.
Now assume \( x_{i} \) and \( y_{i} \) have been constructed.
Take \( t_{i+1} := (x_i + y_i) / 2 \), and compute \( c_n (t_{i+1}) \).
If \( [c_n (u), c_n (t_{i+1})] \cap \Delta = \emptyset \)
then set \( x_{i+1} := t_{i+1} \) and \( y_{i+1} := y_{i} \).
Otherwise, set \( x_{i+1} := x_{i+1} \) and \( y_{i+1} := t_{i+1} \).
It is straightforward to see that the new elements \( x_{i+1} \) and \( y_{i+1} \) also satisfy the properties.
Take the interval \( [u,x_{k}] \) for some relatively large \( k > 0 \), e.g., \( k = 30 \),
for one of the subintervals, say  \( \omega^{\ell} \).
Repeat the same for the other endpoint of \( \omega \) to obtain the other subinterval \( \omega^{r} \), provided that \( c_n (v) \notin \Delta \).

We remark that the convergence of the bisection method is exponential; for example, after \( 30 \) steps, the size of the new interval is computed with the precision of \( 2^{-30} \approx 10^{-9} \) relative to the size of the original interval, which may be satisfactory in most cases. The computation of each step is fast, because it consists of computing \( c_n \) for a single point.
The quantities discussed in the previous two subsections, computed along with the iterations of \( \omega \), can be used further with the smaller intervals, or can be re-computed from scratch; in this paper we chose the second option, because the computation is not very costly, and we can expect to get better estimates for those quantities, due to the smaller interval \( \omega \). Finally, note that due to approximations and rounding, or different monotonicity of \( c_n \) than assumed above, the actual procedure is technically more sophisticated; we discuss the details in Section~\ref{sec:procedure}.

% ====================================================

\section{The Partition}
\label{sec:P}
The construction of the partition \( \mathcal P \), as mentioned above, is based on the computations described in Section \ref{sec:strategies}. However, the way these computations are combined to explicitly construct \( \mathcal P \) is non-trivial, and requires the introduction of some auxiliary constants, and various criteria on when to stop the computations and on how to decide if an interval belongs to \( \mathcal P^{-} \) or \( \mathcal P^{+} \). We describe heuristically, but in some detail,  the overall scheme, and postpone to Section \ref{sec:algorithms} the precise formulation of the formal structure of the algorithms. 

% ====================================================

\subsection{Defining a queue}

The general principle underlying the construction of the partition \( \mathcal P \) is quite simple and is outlined in Section \ref{sec:stocreg}. The construction  relies in a fundamental way on the computations discussed in Section \ref{sec:strategies} and  essentially boils down to a combination of iterating and chopping parameter intervals. We note, however, that this produces a large  number (\emph{possibly millions!}) of small parameter intervals and some criteria need to be put in place regarding the order with which we handle these intervals, at which point we stop iterating, and  how we decide to assign such intervals to either one of the families \( \mathcal P^{+} \) or \( \mathcal P^{-} \). 
For that purpose, we use the notion of a \emph{queue}. In our setting this can be formulated in the following way. At any given moment we have a partition \( \mathcal P \) of \( \Omega \) given by the union of three families of closed intervals with disjoint interiors:
\[
\mathcal P = \mathcal P^{+}\cup \mathcal P^{-}\cup\mathcal P^{q} 
\]
where \( \mathcal P^{+} \) consists of stochastic intervals, \( \mathcal P^{-} \)
 consists of regular intervals, and \( \mathcal P^{q} \) consists of intervals in the queue.  
 Initially, the entire parameter  space \( \Omega \) is placed in the queue as a single interval, and therefore we have 
\begin{equation}\label{eq:Pstart}
\mathcal P = \mathcal P^{q}=\{\Omega\} \quad \text {and } \quad 
 \mathcal P^{+} =  \mathcal P^{-} = \emptyset.
 \end{equation}
 As the process runs, intervals in \( \mathcal P^{q} \) get iterated and  possibly chopped and, according to a set of criteria which we are about to describe, the resulting subintervals are  either assigned to \( \mathcal P^{+} \) or \( \mathcal P^{-} \), after which they are no longer iterated, or to  \( \mathcal P^{q} \) for possible further iteration. Eventually we end up with a situation where 
 \begin{equation}\label{eq:Pend}
  \mathcal P^{q}=\emptyset
  \quad \text {and } \quad 
   \mathcal P = \mathcal P^{+}\cup \mathcal P^{-}
   \end{equation}
at which point we consider to have concluded our construction. In the following subsections we explain the precise mechanism and criteria for moving intervals from the queue into \( \mathcal P^{+} \) or~\( \mathcal P^{-} \) and adding intervals to the queue. We say that an interval is \emph{enqueued} if it is added to the queue, and it is \emph{dequeued} if it is taken back from the queue. 

% ====================================================

\subsection{Processing the queue}
\label{sec:assigning}
The process of moving from the initial partition \eqref{eq:Pstart} to the final partition \eqref{eq:Pend} requires several actions and decisions based on the outcome of the computations described in Section \ref{sec:strategies}. For clarity, we subdivide our explanation of these actions and decisions in a few steps. We suppose that \( \omega\in \mathcal P^{q} \) is an interval in the queue and explain what we do with it and how we decide at some point whether it belongs to \( \mathcal P^{-}  \) or \( \mathcal P^{+} \) or whether it gets chopped, at which point we need to decide what to do with the remaining subintervals. We consider various cases. 

The first and, in some sense, most important case, is when we are successfully able to compute (approximate) iterates of \( \omega \) up to some time \( n\geq 1 \)
 for which \( \omega_{n}\cap \Delta\neq\emptyset \) (or, more precisely, where the outer enclosure of \( \omega_{n} \) intersects \( \Delta \), thus indicating that \( \omega_{n} \) \emph{may}  intersect~\( \Delta \)).   We will consider two subcases.
  
  \medskip\noindent
  \textbf{(P1a)} If \( \omega_{n}\cap \Delta\neq\emptyset \), \( n\geq N_{0} \) and  the escape time conditions \eqref{eq:escape} hold,   we let \( \omega\in \mathcal P^{+} \). 
  
  \medskip\noindent
  This is the one and only situation where we are ``successful'' and place intervals in \( \mathcal P^{+} \). In all other cases below, possibly after subdividing the original interval,  we will either ``give up'' on one or more of the resulting subintervals and place them in \( \mathcal P^{-} \), or save them for further iteration by placing them back in the queue \( \mathcal P^{q} \).

  \medskip\noindent
    \textbf{(P1b)} If \( \omega_{n}\cap \Delta\neq\emptyset \) but \( n< N_{0} \) or  the escape time conditions \eqref{eq:escape}  do not hold, then we chop the interval \( \omega \)  according to the procedure described in Sections \ref{sec:stocreg} and  \ref{sec:avoid} (and in detail in Algorithm~\ref{alg:hitDelta}). 
This chopping procedure subdivides \( \omega \) into at most 3 disjoint subintervals \( \omega=\omega^{\ell}\cup\omega^{\Delta}\cup\omega^{r} \) such that \( \omega^{\ell}_{n}\cap \Delta=\emptyset \) and \( \omega^{r}_{n}\cap \Delta=\emptyset \). We let \( \omega^{\Delta}\in\mathcal P^{-} \) since it intersects~\( \Delta \) and therefore cannot ever satisfy the escape time conditions at any time in the future. The decision about what to do with \( \omega^{\ell}, \omega^{r} \) depends on their size. If they are too small they may contribute little to the final result, and thus one might consider processing them a waste of the computational resources that could otherwise be assigned to investigating larger intervals.   We therefore introduce the variable \( w \geq 0 \) to indicate the minimum width of an interval, relative to the width of \( \Omega \), that we are willing to continue iterating.  If the size of the subintervals is \( \geq w \) we place them back in the queue \( \mathcal P^{q} \), whereas, if their sizes are \( <w \) we ``abandon'' them by placing them in \( \mathcal P^{-} \).  The results described in Section \ref{sec:results} are based on a choice of  \( w=10^{-10} \).

\medskip\noindent
There are only two reasons for which we may fail to arrive at a situation where \( \omega_{n}\cap \Delta\neq\emptyset \): there may be some technical/computational issue which does not allow us to properly compute the iterates of \( \omega \); or it may happen simply that we keep iterating \( \omega \) and it just never hits \( \Delta \). In the first case we distinguish again two subcases.

 \medskip\noindent
  \textbf{(P2a)}
It may happen, possibly due to overestimations caused by the rounding procedures, that for some iterate \( n \) we may have \( 0\in c_{n}'(\omega) \) and/or \( 0\in (f^{n})'(\omega) \) (recall \eqref{eq:cnprime} and \eqref{eq:fnprime}).  The first case indicates a failure of the technical condition \eqref{eq:mon} which is required to continue iterating the interval, and the second a failure of another technical condition which is required to verify some properties of our calculations, see \eqref{eq:ind} in Theorem \ref{thm:intervals} below. 
In both these cases, rather than giving up on \( \omega \) straight away, by placing it in \( \mathcal P^{-} \), we bisect \( \omega \) and consider the two resulting halves of the interval. As in (P1b) we then consisder the size of these subintervals. If they are larger than \( w=10^{-10} \) of the size of \( \Omega \) we place them back in the queue \( P^{q} \), while if they are smaller than \( w \) we place them in \( \mathcal P^{-} \). 

\medskip   \noindent
\textbf{(P2b)} A second technical issue which can arise is the situation where the  lower and upper bounds for the endpoints of \( \omega_n \) are further apart than  the distance between the endpoints themselves. This situation is an effect of rounding errors introduced while evaluating the function~\( f_a \) and suggests that the precision of representable numbers used for the computation is too low. If this happens then we do not get any reasonable lower bound on the width of \( \omega_n \), and therefore iterating \( \omega \) further is pointless; moreover, this situation is explicitly excluded at various steps of our arguments, see \eqref{eq:innerouter} and \eqref{eq:indn}.  There is no way to improve the result, apart from choosing a different precision of numerical computation (choosing a different set \( \mathbf{R} \) of representable numbers). We therefore ``abandon'' such an interval by assigning it to \( \mathcal P^{-} \).  We note, however, that due to the very high precision with which we work, the situation can only occur with extremely small intervals, and therefore this does not seem to provide a significant loss in terms of measure of intervals which eventually make up \( \mathcal P^{+} \).  

\medskip\noindent
Finally, we consider the case where we can continue iterating \( \omega \) but it never intersects \( \Delta \). 

\medskip\noindent
\textbf{(P3)} If \( \omega \) is iterated a huge number of times without hitting \( \Delta \) then this most likely means that the sequence \( \{\omega_n\} \) got trapped inside the attracting neighbourhood of some stable periodic orbit and we are very unlikely to see any escape time in the future. We therefore define the maximum number  \( N_{\max} > 0 \) of iterations that we allow without ever hitting \( \Delta \) and assign an interval to \( \mathcal P^{-} \) if this number if exceeded (see Algorithm~\ref{alg:interval} for the implementation). The results described in Section \ref{sec:results} are based on a choice of  \( N_{\max} = 200 \).

\medskip\noindent
We remark that,  while it is not our goal in this paper to prove that any particular parameters belong to~\( \Omega^{-} \), our rules for placing parameter intervals in \( \mathcal P^{-} \) suggest a strong probability that such intervals belong to, or substantially intersect, open sets in \( \Omega^{-} \). 
 Since we know these intervals explicitly, our calculations may provide the foundations for further work, possibly applying techniques similar to those of~\cite{TucWil09} or~\cite{Gal17}, to actually prove that certain parameter intervals are indeed regular.

%\medskip\noindent
% \textbf{(R2)} Another clue that the interval \( \omega \) may be in a region of parameters in which the corresponding maps have attracting periodic orbits, or in any case that we are unlikely to get any escape times,  may be given by the observation of some subsequence of iterates for which the length of the images \( \omega_{n} \) is decreasing. We implemented a  procedure which compares the length of \( \omega_{n} \) with the length of \( \omega_{n+k} \) when \( \omega_{n} \cap \omega_{n+k} \neq \emptyset \) for some \( n, k > 0 \); if this kind of return of the interval to its previous location is detected a given number \( \chi  \) of times in a row with the same ``period'' \( k \), and the interval gets contracted each time (that is, its length decreases) then the interval \( \omega \) is assigned to \( \mathcal P^{-} \). 
%   The results described in Section \ref{sec:results} are based on a choice of~\( \chi =3 \).

% ====================================================

\subsection{Emptying the queue}
By setting up the numbers \( w \) and \( N_{\max} \), we ensure that all the parameters in \( \Omega \) are eventually moved into either \( \mathcal P^{+} \) or \( \mathcal P^{-} \) and that therefore the process eventually terminates. Our choice of constants used to obtain the results presented in Section \ref{sec:results},  lead to a complete construction of the partition \( \mathcal P \), made up of more than 3.9 million intervals, in only 35 minutes of computation time on a laptop computer. 

It is not completely clear, however, how quickly the computation time may increase if we choose smaller values of the radius \( \delta \) of the critical neighbourhood, larger values of \( N_{0} \), smaller values of the minimum size \( w \) of intervals we consider, or larger values for the number \( N_{\max} \) of iterates before giving up on an interval. It is worth therefore putting in place some ``safeguards'' against the possibility of an essentially never ending computation. We can easily do this by specifying some criteria which limit the amount of computations which we carry out. If any of these criteria are met, we simply stop the computations and transfer all intervals still in the queue to \( \mathcal P^{-} \). This is still completely consistent with the spirit of the construction since the family \( \mathcal P^{-} \) is just the collection of intervals for which we could not verify the escape time condition. The three constraints we can impose are fairly obvious. 

\medskip\noindent
1) We can fix the maximal number \( i_{\max} > 0 \) of intervals to be processed: we keep track of each time an interval is picked form the queue for iteration, until one of the situations described above occurs. After processing this number of intervals, we interrupt the computations.

\medskip\noindent
 2) We can fix the maximal allowed queue size \( q > 0 \); for every interval \( \omega \) that is processed, up to two new intervals are added to the queue when \( \omega_n \) hits \( \Delta \) for some \( n \) or when a problem occurs and the interval \( \omega \) is halved; thus the size of the queue grows linearly during the progress of the computation. If the number of intervals stored in the queue reaches or exceeds \( q \) then we interrupt the computations. Setting the limit on the queue size protects against memory overflow that might be caused by storing too many intervals in the queue, especially if high precision numbers are used that might occupy considerable amount of memory.

\medskip\noindent 
3) We can fix the maximal time \( t > 0 \) (in seconds) that can be used by the program. This constraint is especially useful in order to bound the amount of time that one is willing to wait for the final result, and also to protect the web server's resources when providing access to the program through the web interface.

\medskip

We conclude this section with a discussion of the non-trivial problem of deciding how to prioritize the intervals in the queue, i.e., how to decide which interval to iterate at any given moment. This is especially important if the computation is stopped before the queue is empty, for example, if the program is allowed to run for a limited amount of time only.

From a computational point of view, a  queue is a data structure that is capable of storing objects of certain type, and provides means for extracting them. There are different types of queues in terms of the order in which the objects are extracted. For example, the \emph{fifo queue} (``first in -- first out'') provides the objects in the order in which they were put in the queue (like a typical queue in a supermarket), and the \emph{lifo queue} (``last in -- first out'') provides the most recently stored object first (like a stack of plates).

It seems that a good approach is to use a \emph{priority queue}, in which objects are sorted based on some priority, and the ones with the highest priority are extracted first. More specifically, the queue stores parameter intervals together with the number of times they were successfully iterated, and this number serves as the priority in our queue; \emph{we first extract intervals that were iterated the least number of times}.
In this way, we prioritise those intervals that are lagging behind the others in iterations,
so that we could achieve a state in which all the intervals that are left in the queue
have been iterated at least a certain number of times. if several intervals in the queue have been iterated the same number of times (as is of course often the case) we take the biggest first. In the quite unlikely event that two such intervals are exactly the same size, we introduce other criteria such as priorities depending on the reasons an interval was added to the queue. 

% ====================================================

\subsection{Case study}
\label{sec:study} 

We conclude this section with a case study of the ``history'' of a specific stochastic interval  \( \omega \in \mathcal{P}^+ \) that actually appeared in the computations described in Section~\ref{sec:results}, in order to illustrate some of the processes described above in a concrete case. We consider the interval \( \omega^{(1)} \) for which we have the following outer enclosure when rounding the endpoints to 12 significant digits:
\[
\omega^{(1)} \subset [1.96076793815,1.96077475689].
\]
Its iterates are shown in Table~\ref{tab:bigN}. This is the 1953rd interval taken for iterations from the queue (recall Section~\ref{sec:assigning} for details on how the ``queue'' works). The numbers in Table~\ref{tab:bigN} show that the interval got close to \( \Delta \) at the 7th and 15th iterates. Its width was steadily growing with sudden drops after those two events; eventually the interval ``exploded'' to take up almost the entire phase space at the 26th iterate, thus satisfying the escape time condition with \( |\omega_{N}|\geq 3.5 \).

\begin{table}[htbp]
\centerline{
\begin{tabular}{rcl}
\hline
\( n \) & \( \omega^{(1)}_n \) & \( |\omega^{(1)}_n| \) \\
\hline
0 & \(  [1.9607,1.9608] \) & \( 0.0000068187 \) \\
1 & \( [-1.8839,-1.8838] \) & \( 0.000019921 \) \\
2 & \( [-1.5882,-1.5880] \) & \( 0.000068238 \) \\
3 & \( [-0.56150,-0.56128] \) & \( 0.00020992 \) \\
4 & \( [1.6455,1.6458] \) & \( 0.00022887 \) \\
5 & \( [-0.74766,-0.74689] \) & \( 0.00076011 \) \\
6 & \( [1.4017,1.4030] \) & \( 0.0011428 \) \\
7 & \textcolor{red}{\( [-0.0073968,-0.0041981] \)} & \( 0.0031985 \) \\
8 & \( [1.9607,1.9608] \) & \( 0.000030267 \) \\
9 & \( [-1.8838,-1.8836] \) & \( 0.00012551 \) \\
10 & \( [-1.5879,-1.5873] \) & \( 0.00047967 \) \\
11 & \( [-0.56046,-0.55892] \) & \( 0.0015298 \) \\
12 & \( [1.6466,1.6484] \) & \( 0.0017193 \) \\
13 & \( [-0.75638,-0.75071] \) & \( 0.0056585 \) \\
14 & \( [1.3886,1.3972] \) & \( 0.0085210 \) \\
15 & \textcolor{red}{\( [0.0086113,0.032357] \)} & \( 0.023745 \) \\
16 & \( [1.9597,1.9607] \) & \( 0.00096598 \) \\
17 & \( [-1.8836,-1.8797] \) & \( 0.0037938 \) \\
18 & \( [-1.5871,-1.5727] \) & \( 0.014284 \) \\
19 & \( [-0.55781,-0.51266] \) & \( 0.045141 \) \\
20 & \( [1.6496,1.6980] \) & \( 0.048329 \) \\
21 & \( [-0.92227,-0.76048] \) & \( 0.16177 \) \\
22 & \( [1.1101,1.3825] \) & \( 0.27222 \) \\
23 & \( [0.049666,0.72824] \) & \( 0.67856 \) \\
24 & \( [1.4304,1.9584] \) & \( 0.52784 \) \\
25 & \( [-1.8742,-0.085416] \) & \( 1.7887 \) \\
26 & \( [-1.5518,1.9535] \) & \( 3.5052 \) \\
\hline
\end{tabular}
}
\caption{Iterates of one specific interval \( \omega  \in \mathcal{P}^+ \) with large \( |\omega_n| \) at escape time.
All the numbers rounded to 5 significant digits. An outer bound on each \( \omega_n \) is shown,
as well as a lower bound on its width, as calculated in the high-precision arithmetic.
Close encounters with \( \Delta \) are shown in red. Full discussion in Section \ref{sec:study}.}
\label{tab:bigN}
\end{table}

Let us check the circumstances under which this interval entered the queue. Each interval that is put in the queue is assigned a consecutive number, starting from 1 that was assigned to the original interval \( \Omega = [1.4,2] \). The computation log shows that \( \omega^{(1)} \) was assigned the number 2565, and it was put in the queue as a result of halving another interval, let us call it its parent and denote by \( \omega^{(0)} \). This was the 1914th processed interval, and it was halved due to a problem with determining the sign of \( c'_n(\omega) \), discussed in Section~\ref{sec:assigning} as subcase (P2a), after having computed its 7th iterate. Recall that \( \omega^{(1)}_7 \) was indeed close to \( \Delta \); which means that halving the interval \( \omega^{(0)} \) instead of throwing it into \( \mathcal{P}^- \) was a good decision, because it allowed saving at least a half of it for \( \mathcal{P}^+ \).

Let us have a look at the ``twin brother'' \( \omega^{(2)} \) of the successful interval \( \omega^{(1)} \). Its 12-digit outer enclosure is \( [1.96076111942,1.96076793816] \); it was assigned the number 2564 when it was put in the queue. It was iterated just before \( \omega^{(1)} \), that is, it was the 1952nd iterated interval. It was subject to the same problem as \( \omega^{(0)} \) and was halved after 7 iterates. Its two halves \( \omega^{(3)} \) and \( \omega^{(4)} \) were put back in the queue, got consecutive numbers 2619 and 2620, and were iterated as the 1992nd and 1993rd intervals, respectively. The interval \( \omega^{(4)} \) hit \( \Delta \) after a remarkable number of 31 iterates, and the width of its 31st iterate slightly exceeded \( 1.5 \), so it was added to \( \mathcal{P}^+ \). This means that we already qualified \( 3/4 \) of \( \omega^{(0)} \) as stochastic! However, the problem persisted for \( \omega^{(3)} \), which was then halved again. Its children \( \omega^{(5)} \) and \( \omega^{(6)} \) got consecutive numbers 2678 and 2679 in the queue. They were pulled out from the queue and processed as the 2032nd and 2033rd intervals, respectively. The problem persisted for \( \omega^{(5)} \), while \( \omega^{(6)} \) was iterated 15 times until it hit \( \Delta \); its 15th iterate was contained in \( [0.000554,0.002311] \) and was not large enough to put the interval in \( \mathcal{P}^+ \), so it was chopped; note that the left endpoint of \( \omega^{(6)}_{15} \) was actually in \( \Delta \), so the chopping resulted in only one part put back in the queue. The loss was considerable, only \( 70\% \) of \( \omega^{(6)} \) survived the chopping.
We stop our investigation here. Although the fraction of \( \omega^{(0)} \) qualified as stochastic did not increase to \( 7/8 \), there is hope that some of the descendants of \( \omega^{(6)} \) eventually contributed to \( \mathcal{P}^+ \) in further iterations.

This short excerpt of the family saga of the interval \( \omega^{(1)} \) illustrates the main ideas of our approach in constructing the sets \( \mathcal{P}^+ \) and \( \mathcal{P}^- \), and shows a variety of dynamical situations encountered.

% ====================================================

\section{The Numerics}
\label{sec:procedure}

The strategies introduced in Section~\ref{sec:strategies} are intertwined into a single computational procedure for computing inner and outer bounds for \( \omega_n \) together with rigorous estimates for the derivatives \eqref{eq:chain} and \eqref{eq:quot}, and splitting the interval \( \omega \) into smaller parts whenever condition \eqref{eq:mon} is not satisfied or \( \omega_n \) hits \( \Delta \).
In this section, we explain the issues involved in making sure we obtain \emph{rigorous} bounds for all these calculations. In Section \ref{sec:algorithms} we then describe the structure of the algorithms used to implement the calculations. 

% ====================================================

\subsection{Precision}

The very first step in the construction is the choice of the set of representable real numbers \( \mathbf{R} \subset \mathbb{R} \).
In practice, the choice of this set depends on the representation of numbers used in the software,
and is different for double-precision floating point numbers following the IEEE~754 standard \cite{ieee754}
than for floating-point numbers of fixed size implemented by the GNU MPFR software library \cite{mpfr}.
For clarity of presentation, within Section~\ref{sec:algorithms}
we are going to use bold typeface to denote elements of \( \mathbf{R} \);
for example, \( \mathbf{a} \in \mathbf{R} \), as opposed to the general \( a \in \mathbb{R} \).

It is important to be aware of the fact that the actual result of an arithmetic operation or the result of the computation of the value of a function on representable numbers need not be a representable number in general. However, in a proper setting, it is possible to request that such results are rounded downwards or upwards to representable numbers in the actual machine computations. Therefore, even if it is not possible in general to compute the exact value of many expressions, it is always possible to compute a lower and an upper bound for each of them.
In the case of elementary operations, such as addition or multiplication, the standards such as IEEE~754 typically require that the result is rounded downwards or upwards to the closest representable number in the corresponding direction; thanks to this feature, the inaccuracy caused by rounding is minimised.

The most important quantity of interest regarding the precision of our calculations is  the binary precision \( p  \) to be used for representable numbers implemented by the MPFR library; for example, if \( p = 250 \) then the relative accuracy of numbers used in the computations is roughly \( 2^{-250} \approx 10^{-80} \), that is, all the numbers are rounded at the \( 80 \)th decimal digit, which we consider quite high precision for the results, yet computationally feasible in terms of the speed and memory usage.  The choice of binary precision is closely related to the number of iterations which we want to compute, and \( p=250 \) is quite sufficient for the number of iterations we consider for the results we describe in Section \ref{sec:results}. Higher precision would be desirable and could easily be implemented if we carried out the calculations for higher values of \( N_{0} \). 

A second quantity which  affects to some extent the precision of our calculations is the number \( s  \) of bisection steps used in the chopping procedure described in Section \ref{sec:avoid} (and more formally in  Algorithm~\ref{alg:bisection} below); the higher the value, the more expensive the computations; the lower the value, the less accurate the chopping procedure. In order to determine reasonable balance between these two, we made some experiments to check the improvement in the results and the increase of computation time with the increase in \( s \), and we decided to use \( s = 40 \).

% ====================================================

\subsection{Rounding}
\label{sec:repr}
Instead of introducing separate notation for representable versions of all the operations and functions
with rounding downwards or upwards, we use the two special assignment symbols ``\( :\leq \)'' and ``\( :\geq \)''
instead of ``\( := \)'' in order to indicate the rounding direction.
Specifically, if \( \varphi \colon \mathbb{R}^k \to \mathbb{R} \)
for some \( k \in \mathbb{N} \) then
\[
\mathbf{u} :\leq \varphi (\mathbf{x}_1, \ldots, \mathbf{x}_k)
\]
means that \( \mathbf{u} \in \mathbf{R} \) is the result of machine computation of a representable number that is a lower bound for the actual value of \( \varphi (\mathbf{x}_1, \ldots, \mathbf{x}_k) \).
We define the upwards-rounded counterpart
\[
\mathbf{v} :\geq \varphi (\mathbf{x}_1, \ldots, \mathbf{x}_k)
\]
in an analogous way.
If the direction of rounding is not important, we use the ``rounding to the nearest'' mode, and we use the symbol \( :\approx \) to explicitly indicate the fact that rounding to a representable number takes place when computing the expression on the right hand side of \( :\approx \). This happens, for example, when we compute an approximation of the middle of an interval: ``\( \mathbf{c} :\approx (\mathbf{a} + \mathbf{b}) / 2 \).''
We remark that even if one computes a lower bound, the rounding direction does not have to be ``downwards'' in all the operations. Consider, for example, the computation of \( 1 / (x + y) \) for \( x, y > 0 \). One would first compute \( \mathbf{z} :\geq x + y \) and then \( \mathbf{s} :\leq 1 / \mathbf{z} \) so that the number \( \mathbf{s} \) is a lower bound on \( 1 / (x + y) \).

All the numbers that appear in the algorithms and in the computations must be representable. In case any number appears in the description that is not exactly reprsentable in the binary floating-point arithmetic that we use, such as \( 1.4 \) for example, it is implicitly rounded to the nearest representable number when passed to the algorithm.

\medskip

Before moving on to describe our main procedure, we introduce some notation. 
\begin{definition}\label{def:defsign}
An interval \( [x^-,x^+] \) is said to be \emph{of definite sign} if \( 0 \notin [x^-,x^+] \).
\end{definition}
Note that an interval is of definite sign if its both endpoints are not zero and have the same sign.
Given two intervals \( x = [x^-,x^+] \) and \( y = [y^-,y^+] \),
let 
\[
g^- (x,y):=\min\{-2uv : u \in x, v \in y \} \quad\text{ and } \quad 
g^+ (x,y):=\max\{-2uv : u \in x, v \in y \}
\]
These functions provide the tightest outer enclosure for the arithmetic operation \( -2xy \) on intervals. Notice that the intervals \( x\) and \( y  \) are not necessarily representable and even if they were, the output of the functions \( g^{-}, g^{+} \) are not necessarily representable. However, if they are of definite sign, they are given by a simple formula 
\begin{equation*}
\label{eq:g}
g^-(x,y) =
\begin{cases}
-2 x^+ y^+ & \text{if } x > 0, y > 0, \\
-2 x^- y^- & \text{if } x < 0, y < 0, \\
-2 x^- y^+ & \text{if } x > 0, y < 0, \\
-2 x^+ y^- & \text{if } x < 0, y > 0;
\end{cases}
\qquad
g^+(x,y) =
\begin{cases}
-2 x^- y^- & \text{if } x > 0, y > 0, \\
-2 x^+ y^+ & \text{if } x < 0, y < 0, \\
-2 x^+ y^- & \text{if } x > 0, y < 0, \\
-2 x^- y^+ & \text{if } x < 0, y > 0.
\end{cases}
\end{equation*}
and can be rounded up and down to get inner and outer enclosures of the set \( \{-2uv : u \in x, v \in y \} \) by representable  numbers. 

% ====================================================

\subsection{Inductive assumptions}
\label{sec:indStart}

In this subsection, we introduce an inductive procedure for iterating \( \omega \) and computing a sequence of numbers that, as we shall prove,  provide lower and upper bounds for the quantities discussed in Section~\ref{sec:strategies}.
Let \( \mathbf{a} < \mathbf{b} \) be representable numbers, and consider the parameter interval
\begin{equation*}
\omega := [\mathbf{a},\mathbf{b}].
\end{equation*}
We are going to construct  bounds on the quantities discussed in Sections \ref{sec:iterate} and~\ref{sec:estimate}. To formulate these bounds we will define inductively representable intervals
\begin{equation}\label{eq:int}
\mathbf{a}_n := [\mathbf{a}_n^-,\mathbf{a}_n^+], \ 
\mathbf{b}_n := [\mathbf{b}_n^-,\mathbf{b}_n^+],  \ 
\mathbf{c}_n := [\mathbf{c}_n^-,\mathbf{c}_n^+], \
\mathbf{d}_n := [\mathbf{d}_n^-,\mathbf{d}_n^+], \
\mathbf{f}_n := [\mathbf{f}_n^-,\mathbf{f}_n^+]
\end{equation}
and (when possible) two additional intervals
\begin{equation}\label{eq:omegadef}
\underline{\omega_{n}}\subseteq \overline{\omega_{n}},
\end{equation}
defined in terms of those above, where \( \overline{\omega_{n}} \) is the \emph{convex hull} of \( \mathbf{a}_n \) and \( \mathbf{b}_n \), i.e., the smallest closed interval containing both \( \mathbf{a}_n \) and \( \mathbf{b}_n \), and \( \underline{\omega_{n}} \) is the closure of the \emph{unique bounded component} of \( \mathbb R\setminus \{\mathbf{a}_n\cup \mathbf{b}_n\} \) \emph{if} \( \mathbf{a}_n\cap \mathbf{b}_n=\emptyset  \) (and \emph{undefined} otherwise). Clearly \( \underline{\omega_{n}}, \overline{\omega_{n}} \), when defined, are also representable intervals, and if \( \mathbf{a}_n\cap \mathbf{b}_n=\emptyset  \), are given by the following simple formulas: 
\begin{equation}
\label{eq:innerouter}
\underline{\omega_n} :=
\begin{cases}
[\mathbf{a}_n^+, \mathbf{b}_n^-] & \text{if } \mathbf{a}_n^+ < \mathbf{b}_n^-, \\
[\mathbf{b}_n^+, \mathbf{a}_n^-] & \text{otherwise}.
\end{cases}
\qquad \text{ and } \qquad
\overline{\omega_n} :=
\begin{cases}
[\mathbf{a}_n^-, \mathbf{b}_n^+] & \text{if } \mathbf{a}_n^- < \mathbf{b}_n^+, \\
[\mathbf{b}_n^-, \mathbf{a}_n^+] & \text{otherwise}.
\end{cases}
\end{equation}
The definition of the intervals \eqref{eq:int} is inductive and 
for the initialisation of the induction, \( n=0 \), we define the following representable numbers:
\begin{equation}
\label{eq:ind0}
\mathbf{a}_0^- \coloneqq \mathbf{a}_0^+ \coloneqq \mathbf{a},\quad
\mathbf{b}_0^- \coloneqq \mathbf{b}_0^+ \coloneqq \mathbf{b},\quad
\mathbf{c}_{0}^-= \mathbf{c}_{0}^+ \coloneqq 1, \quad
\mathbf{f}_{0}^-= \mathbf{f}_{0}^+ \coloneqq 1,\quad
\mathbf{d}_{0}^-= \mathbf{d}_{0}^+ \coloneqq 1,
\end{equation}
and define the corresponding intervals \( \mathbf{a}_0, \mathbf{b}_0, \mathbf{c}_0,\mathbf{d}_0,\mathbf{f}_0, \overline{\omega_{0}}, \underline{\omega_{0}}  \) as in \eqref{eq:int}-\eqref{eq:omegadef}. 
Note that we admit degenerate intervals as singletons if both endpoints are equal, and we distinguish such intervals (sets) from the individual numbers such as \( \mathbf{a} \) and \( \mathbf{b} \). We also consider our intervals ``ordered'' in the sense that the left endpoint as written is always assumed to be \( \leq \) the right endpoint. Let us now assume inductively that intervals \( \mathbf{a}_k, \mathbf{b}_k, \mathbf{c}_k,\mathbf{d}_k,\mathbf{f}_k\), and therefore also  \(  \overline{\omega_{k}} \) (but not necessarily \( \underline{\omega_{k}}  \)), have been defined for some \( k\geq 0 \), \emph{and that \( \overline{\omega_{k}} \) has definite sign}: 
\begin{equation}\label{eq:defsign}\tag*{\ensuremath{(\star)_{k}}}
0\notin\overline{\omega_{k}}.
\end{equation}
With these assumptions, we define the intervals \( \mathbf{a}_{k+1}, \mathbf{b}_{k+1}, \mathbf{c}_{k+1},\mathbf{d}_{k+1},\mathbf{f}_{k+1}  \) in the next section.

% ====================================================

\subsection{Inductive step}
\label{sec:indStep}

Recall that \( \mathbf{a} \) and \( \mathbf{b} \) (without subscripts) are the endpoints of the interval \( \omega \) and for every \( a\in \omega \),  \( f_a \) is the quadratic map defined in \eqref{eq:quadratic}.
Recall from Section~\ref{sec:repr} that we use the notation ``\( \mathbf{u} :\leq \varphi(\mathbf{y}) \)'' to indicate that \( \mathbf{u} \) is computed to be a lower bound on the expression \( \varphi(\mathbf{y}) \); similarly with ``\( :\geq \)''.
If \( \mathbf{a}_k^+ < 0 \), we set
\begin{align}
\label{eq:indIncreasing}
&\mathbf{a}_{k+1}^- :\leq f_{\mathbf{a}}(\mathbf{a}_{k}^-),\quad
\mathbf{a}_{k+1}^+ :\geq f_{\mathbf{a}}(\mathbf{a}_{k}^+),\quad
\mathbf{b}_{k+1}^- :\leq f_{\mathbf{b}}(\mathbf{b}_{k}^-),\quad
\mathbf{b}_{k+1}^+ :\geq f_{\mathbf{b}}(\mathbf{b}_{k}^+);
\end{align}
otherwise, we set % If \( \mathbf{a}_n^- \ge 0 \), we set
\begin{align}
\label{eq:indDecreasing}
&\mathbf{a}_{k+1}^- :\leq f_{\mathbf{a}}(\mathbf{a}_{k}^+),\quad
\mathbf{a}_{k+1}^+ :\geq f_{\mathbf{a}}(\mathbf{a}_{k}^-),\quad
\mathbf{b}_{k+1}^- :\leq f_{\mathbf{b}}(\mathbf{b}_{k}^+),\quad
\mathbf{b}_{k+1}^+ :\geq f_{\mathbf{b}}(\mathbf{b}_{k}^-).
\end{align}
Then we let
\begin{align}
%&\mathbf{a}_n\coloneqq \min\{c_n(\mathbf{a}),c_n(\mathbf{b})\},\quad \mathbf{b}_n\coloneqq \max\{c_n(\mathbf{a}),c_n(\mathbf{b})\},\\
\label{eq:indC}
\mathbf{c}_{k+1}^- &:\leq 1 + g^- (\mathbf{c}_k,\ovl{\omega_k}),
\qquad
\mathbf{c}_{k+1}^+ :\geq 1+ g^+(\mathbf{c}_k,\ovl{\omega_k}), \\
\label{eq:indF}
 \mathbf{f}_{k+1}^- &:\leq g^-(\mathbf{f}_{k}, \ovl{\omega_k}),
 \qquad \qquad
\mathbf{f}_{k+1}^+ :\geq g^+(\mathbf{f}_{k}, \ovl{\omega_k} ),\\
\label{eq:indD}
\mathbf{d}_{k+1}^- &:\leq \mathbf{d}_{k}^- + 1/{\mathbf{f}_{k+1}^+}, 
\qquad\quad
\mathbf{d}_{k+1}^+ :\geq \mathbf{d}_{k}^+ + 1/{\mathbf{f}_{k+1}^-}.
\end{align}
It is easy to see that this gives well defined intervals \( \mathbf{a}_{k+1}, \mathbf{b}_{k+1}, \mathbf{c}_{k+1},\mathbf{d}_{k+1},\mathbf{f}_{k+1}  \) and that these are explicitly and rigorously computable given the representable intervals \( \mathbf{a}_k, \mathbf{b}_k, \mathbf{c}_k,\mathbf{d}_k,\mathbf{f}_k\) and under assumption \( (\star)_{k} \). In fact, \( (\star)_{k} \) is only required to ensure that the intervals \( \mathbf{a}_{k+1}, \mathbf{b}_{k+1} \) defined in \eqref{eq:indIncreasing}-\eqref{eq:indDecreasing} are well-defined, the other intervals are well-defined with no assumptions. 
However, it is not immediate, nor in fact is it always the case, that these intervals give us any dynamical information. In the next section we prove a non-trivial result which gives conditions for these intervals to provide the required bounds. 

% ====================================================

\subsection{Rigorous bounds}
The main result of this section gives conditions which ensure that the intervals defined above  provide the bounds for the required quantities.

\begin{theorem}
\label{thm:intervals}
Let \( \omega := [\mathbf{a},\mathbf{b}]\) be a parameter interval and let  \( n \geq 0 \). Suppose 
 that for every \( 0\leq k \leq n \) the intervals  \( \mathbf{a}_k, \mathbf{b}_k, \mathbf{c}_k,\mathbf{d}_k,\mathbf{f}_k \) have been defined as above and satisfy condition \( (\star)_{k} \). Then \( \mathbf{a}_{n+1}, \mathbf{b}_{n+1}, \mathbf{c}_{n+1},\mathbf{d}_{n+1},\mathbf{f}_{n+1}  \) are defined and
\begin{equation}\label{eq:indab}
 c_{n+1} (\mathbf{a})\in \mathbf{a}_{n+1} 
 \quad\text{ and } \quad 
c_{n+1} (\mathbf{b}) \in \mathbf{b}_{n+1}. 
\end{equation}
If \( \mathbf{c}_{n}\) and \( \mathbf{f}_{n} \) have definite sign, then also 
\begin{equation}\label{eq:ind}
 \quad
c_{n+1}' (\omega) \subseteq \mathbf{c}_{n+1},  \quad 
(f^{n+1})' (\omega) \subseteq \mathbf{f}_{n+1},\quad 
 {c_{n+1}'}/{(f^{n+1})'} \subseteq \mathbf{d}_{n+1}. 
\end{equation}
If, moreover,   \( \mathbf{a}_{n+1}\cap \mathbf{b}_{n+1} =\emptyset  \)  
and \(  \mathbf{c}_{n+1} \) has definite sign, then 
\begin{equation}\label{eq:indn}
 \underline{{\omega}_{n+1}}\subseteq {\omega}_{n+1}\subseteq \overline{{\omega}_{n+1}}.
\end{equation}
\end{theorem}

We emphasise the fact that the assumptions of the theorem for a fixed \( n \geq 0 \) can be verified by means of finite machine computation, which includes the computation of the various numbers and intervals. In the next section, we introduce Algorithm~\ref{alg:interval} that does precisely this. The conclusion of the theorem, however, provides nontrivial mathematical properties whose verification may not be obvious at all. In particular, the inner and outer bounds on \( \omega_n \), and the outer bounds on the various derivatives, computed in the inductive way using the formulas provided in this subsection, are nontrivial ingredients of the computations needed in the bigger project outlined in Section~\ref{sec:motivation}.

\begin{proof}
We prove Theorem~\ref{thm:intervals} by induction on  \( n \).
For \( n = 0 \), \eqref{eq:indab} and \eqref{eq:ind} follow directly from \eqref{eq:ind0}
and the fact that for every \( a\in\omega \), we have \( c_0(a)=a \), and thus \( c_0'(a)=1 \), and 
moreover, \( f^0 \) is the identity map, so \( (f^0)'(a) = 1 \).
Condition \eqref{eq:indn} follows immediately from \eqref{eq:ind0} and \eqref{eq:innerouter},
and from the fact that \( \mathbf{a} < \mathbf{b} \). 
We therefore assume inductively the conclusions of the theorem for all \( 0\leq k \leq n \) under the corresponding assumptions. 

To prove \eqref{eq:indab}, since \( \overline{\omega_n} \) is of definite sign, \( f_a \) is monotone on \( \overline{\omega_k} \) for every \( a \), in particular for \(a = \mathbf{a} \) and for \(a = \mathbf{b} \). Since \( \mathbf{a}_n \subset \overline{\omega_n} \), the direction of this monotonicity can be determined by the single number \( \mathbf{a}_n^+ \). If this number is negative then both \( f_{\mathbf{a}} \) and \( f_{\mathbf{b}} \) are increasing on \( \overline{\omega_n} \), and the numbers computed using \eqref{eq:indIncreasing} satisfy \( \mathbf{a}_{n+1}^- \leq \mathbf{a}_{n+1}^+ \) and \( \mathbf{b}_{n+1}^- \leq \mathbf{b}_{n+1}^+ \); the reasoning is analogous if \eqref{eq:indDecreasing} has to be used.
This argument, combined with the formula \( c_{n+1} (\mathbf{a}) = f_{\mathbf{a}} (c_n (\mathbf{a})) \) that defines the critical orbit for \( \mathbf{a} \) (and similarly for \( \mathbf{b} \)), proves \(  c_{n+1} (\mathbf{a})\in \mathbf{a}_{n+1} \)  and 
\( c_{n+1} (\mathbf{b}) \in \mathbf{b}_{n+1}\).

The three terms in \eqref{eq:ind}  are all proved by almost the same argument. For the first one, 
notice that the formula \eqref{eq:indC} can be seen as a numerical version of \eqref{eq:derivC};
specifically, if \( \mathbf{c}_{n} \) is an interval containing \( c'_{n}(a) \), as per our inductive assumptions, and if \( \mathbf c_{n} \) and \( \overline{\omega_{n}} \) are of definite sign, as per the assumptions in the theorem, then \eqref{eq:indC} provides an interval that contains \( c'_{n+1}(a) \), i.e. \(c_{n+1}' (\omega) \subseteq \mathbf{c}_{n+1}\). A very similar argument applies to the last two terms of \eqref{eq:ind}  except we look at \eqref{eq:indF} as a numerical version of \eqref{eq:chain}, and \eqref{eq:indD} as a numerical version of \eqref{eq:quot}.

Finally,  to prove \eqref{eq:indn}, first notice that \( c_{n+1} \) is monotone on \( \omega \): this follows from the fact that \( \mathbf{c}_{n+1} \) is of definite sign, as per the assumptions of the theorem,  combined with the just proved property \eqref{eq:ind} stating that \( c'_{n+1}(\omega) \subseteq \mathbf{c}_{n+1} \). Thanks to this monotonicity, the image of the interval \( \omega \) by \( c_{n+1} \) lies entirely between the images of its endpoints, \( \mathbf{a} \) and \( \mathbf{b} \). The images of these points are contained in the corresponding intervals \( \mathbf{a}_{n+1} \) and \( \mathbf{b}_{n+1} \), respectively; the latter fact was just proved as \eqref{eq:indab}. Under the assumption that \( \mathbf{a}_{n+1} \) and \( \mathbf{b}_{n+1} \) are disjoint  the formula \eqref{eq:innerouter} clearly defines \( \overline{\omega_n} \) as the smallest interval containing both intervals \( \mathbf{a} \) and \( \mathbf{b} \), and therefore containing both endpoints of \( c_{n+1} (\omega) \), and thus the entire interval \( c_{n+1} (\omega) \). Moreover, \eqref{eq:innerouter} defines \( \underline{\omega_{n+1}} \) as the closure of \( \overline{\omega_{n+1}} \setminus (\mathbf{a} \cup \mathbf{b}) \), which is an interval contained between the images of the endpoints of \( \omega \), and therefore contained in \( c_{n+1} (\omega) \).
\end{proof}

% ====================================================

\section{The Algorithms}
\label{sec:algorithms}

% formatting!
\vskip -3pt

In this section, we introduce algorithms that serve the purpose of conducting the computations described in Section~\ref{sec:procedure}.
While introducing the algorithms, we are going to use the concept of a \emph{controller}. It is an object to which the progress of computations is reported, which submits obtained results for further processing if desired, and which is responsible for making decisions on how to proceed whenever problems are encountered.

% ====================================================

\subsection{Algorithm for iterating a parameter interval}
\label{sec:algorithm}

% formatting!
\vskip -3pt

Algorithm~\ref{alg:interval} below conducts inductive computations described in Sections \ref{sec:indStart}--\ref{sec:indStep} for a single interval \( \omega \) of paramters, and verifies the assumptions of Theorem~\ref{thm:intervals} at each iteration. The algorithm is defined in the form of an iterative procedure that is in principle indefinite; therefore, in the actual computations described in Section~\ref{sec:results}, we impose some specific stopping criteria that are enforced by the controller. Note that there is no single object returned by the algorithm as its output; instead, the algorithm produces a multitude of data, and supplies this data to the controller that might, for example, store it in a file or send to another procedure for further processing. The details on how the controller reacts to the different events indicated by calling its various functions in Algorithm~\ref{alg:interval} are discussed and explained in Section~\ref{sec:problems}, and the procedure for splitting the interval \( \omega \) if its iteration hits the critical neighbourhood \( \Delta = (-\delta,\delta) \) is provided in Section~\ref{sec:chop}. The instruction ``break'' makes the algorithm exit the loop.

% Nonindent is necessary here for hyperref to produce a correct hyper link, see explanation here: https://tex.stackexchange.com/questions/14000/enumerate-after-label
\begin{algorithm}\noindent\label{alg:interval}
\rm
\begin{tabbing}
\hspace{.5 cm}\= \hspace{.5 cm}\= \hspace{.5 cm}\= \hspace{.5 cm}\= %
\hspace{.5 cm}\= \hspace{.5 cm}\= \hspace{.5 cm}\= \hspace{.5 cm}\= \kill
\kw{function} process\_an\_interval \\
\kw{input:} \\
\> \( \omega = [\mathbf{a}, \mathbf{b}] \): an interval; \\
\kw{begin} \\
\> initialize the induction as defined by \eqref{eq:ind0}; \\
\> define \( \overline{\omega_0} \) and \( \underline{\omega_0} \) following \eqref{eq:innerouter}; \\
\> \kw{for} \( n := 0, 1, 2, 3, \ldots \) \kw{do}: \\
\>\> compute \( \mathbf{c}_{n+1} \) following \eqref{eq:indC}; \\
\>\> \kw{if} \( 0 \in \mathbf{c}_{n+1} \) \kw{then} \\
\>\>\> controller.problemC (\( \omega \), \( n \)); \kw{break}; \\
\>\> compute \( \mathbf{a}_{n+1} \) and \( \mathbf{b}_{n+1} \) following \eqref{eq:indIncreasing} or \eqref{eq:indDecreasing}, as appropriate; \\
\>\> \kw{if} \( \mathbf{a}_{n+1} \cap \mathbf{b}_{n+1} \neq \emptyset \) \kw{then} \\
\>\>\> controller.innerEmpty (\( \omega \), \( n \)); \kw{break}; \\
\>\> define \( \overline{\omega_{n+1}} \) and \( \underline{\omega_{n+1}} \) following \eqref{eq:innerouter}; \\
\>\> compute \( \mathbf{f}_{n+1} \) following \eqref{eq:indF}; \\
\>\> \kw{if} \( 0 \in \mathbf{f}_{n+1} \) \kw{then} \\
\>\>\> controller.problemF (\( \omega \), \( n \)); \kw{break}; \\
\>\> compute \( \mathbf{d}_{n+1} \) following \eqref{eq:indD}; \\
\>\> controller.notify (\(\omega\), \( n+1 \)); \\
\>\> \kw{if} \( \inter \overline{\omega_{n+1}} \cap \Delta \neq \emptyset \) \kw{then} \\
\>\>\> omegaHitDelta (\( \omega \), \( n+1 \)); \kw{break}; \\
\kw{end.}
\end{tabbing}
\end{algorithm}

The next result is an immediate consequence of the fact that Algorithm~\ref{alg:interval} follows the construction introduced in Sections \ref{sec:indStart}--\ref{sec:indStep} and verifies the assumptions of Theorem~\ref{thm:intervals}, except instead of checking that the interval \( \overline{\omega_{n+1}} \) is of definite sign, it checks a stronger condition; namely, given some \( \delta > 0 \), the algorithm verifies whether the distance of \( \overline{\omega_{n+1}} \) from the critical point \( c = 0 \) is at least \( \delta \).

\begin{corollary}
\label{cor:intervals}
Let Algorithm \ref{alg:interval} be called with a compact interval \( \omega = [\mathbf{a},\mathbf{b}] \subseteq [1,\infty) \) with \( \mathbf{a} < \mathbf{b} \). Assume that the radius \( \delta > 0 \) of the critical neighbourhood satisfies \( \delta < 1 \).
Then, every time the algorithm makes a call to the procedure \emph{controller.notify} with \( \omega \) and \( n + 1 \),
the quantities computed by this algorithm satisfy the properties \eqref{eq:indab}--\eqref{eq:indn}.
\end{corollary}

% ====================================================

\subsection{A queue of parameter intervals}
\label{sec:queue}

In this section we introduce an algorithm for managing a collection of intervals that are waiting to be processed using Algorithm~\ref{alg:interval}.
The intervals in some input collection are first added to the queue, obviously together with the number of times they were previously iterated defined as~\( 0 \).
We assume that the controller introduced in Algorithm~\ref{alg:interval} has unlimited access to this queue.
In the framework of the computations, intervals are extracted from the queue one by one and processed individually by Algorithm~\ref{alg:interval}.
This procedure is introduced in Algorithm \ref{alg:queue} below.

It is important to mention here that some intervals with certain priorities might be added to the queue by the controller, in response to the different situations that may be encountered in Algorithm~\ref{alg:interval}.
This feature makes the problem of determining which interval to process more sophisticated than just sorting the list of initial intervals at the beginning and processing them in this order. This observation justifies using the structure of a queue for that purpose.

% Nonindent is necessary here for hyperref to produce a correct hyper link, see explanation here: https://tex.stackexchange.com/questions/14000/enumerate-after-label
\begin{algorithm}\noindent\label{alg:queue}
\rm
\begin{tabbing}
\hspace{.5 cm}\= \hspace{.5 cm}\= \hspace{.5 cm}\= \hspace{.5 cm}\= %
\hspace{.5 cm}\= \hspace{.5 cm}\= \hspace{.5 cm}\= \hspace{.5 cm}\= \kill
\kw{function} process\_all\_intervals \\
\kw{input:} \\
\> \( \{\omega^i = [\mathbf{a}^i, \mathbf{b}^i]\}_{i = 1}^{M} \) for some natural \( M \ge 0 \) \\
\kw{begin} \\
\> \( Q \) := a queue of pairs (interval, integer); \\
\> \kw{for} \( i := 1 \) to \( M \): \\
\>\> \( Q \).enqueue \( (\omega^i \), \( 0 \)); \\
\> \kw{while} \( Q \) is not empty: \\
\>\> \( \omega \) := \( Q \).dequeue(); \\
\>\> process\_an\_interval (\( \omega \)); \emph{// Algorithm~\ref{alg:interval}} \\
\kw{end.}
\end{tabbing}
\end{algorithm}

% ====================================================

\subsection{Overestimate problems}
\label{sec:problems}

Whenever assumptions of Theorem~\ref{thm:intervals} cannot be successfully verified in Algorithm~\ref{alg:interval}, the controller is notified and must take some action. The two obvious choices are either to abandon the problematic interval and not to consider it for further processing, or subdivide it into smaller parts and put some or all of them in the queue \( Q \) that is defined in Algorithm~\ref{alg:queue} above. In this section we describe the actions that we chose to undertake in the cases shown in Algorithm~\ref{alg:interval}.

The two problems with verifying the various technical assumptions in Algorithm~\ref{alg:interval}, reported to the controller using the functions \emph{controller.problemC} and \emph{controller.problemF}, are of similar nature. The first problem is reported when we fail to verify \eqref{eq:mon} that would imply the monotonicity of \( c_n \), and thus we cannot use our method for computing a rigorous bound on \( \omega_n \) introduced in Section~\ref{sec:iterate}. It is likely that the problem with verifying the monotonicity of \( c_n \) is caused in many cases by considerable overestimates in computing the rigorous bound \( \mathbf{c}_n \) for \( c'_n (\omega) \). The second problem appears if the overestimates in computing an outer bound \( \mathbf{f}_n \) for \( (f^{n})'(\omega) \) become so bad that the bound includes \( 0 \), which is obviously wrong. Algorithm~\ref{alg:halve} shows a suggestion of what one can do in these two situations. Our strategy is to halve the interval \( \omega \) in hope that the problem will disappear (which indeed often happens, as illustrated in the case study described in Section~\ref{sec:study}). The controller puts both halves of \( \omega \) to the queue \( Q \) so that these smaller intervals can be processed later.

% Nonindent is necessary here for hyperref to produce a correct hyper link, see explanation here: https://tex.stackexchange.com/questions/14000/enumerate-after-label
\begin{algorithm}\noindent\label{alg:halve}
\rm
\begin{tabbing}
\hspace{.5 cm}\= \hspace{.5 cm}\= \hspace{.5 cm}\= \hspace{.5 cm}\= %
\hspace{.5 cm}\= \hspace{.5 cm}\= \hspace{.5 cm}\= \hspace{.5 cm}\= \kill
\kw{function} controller.problemC, \kw{function} controller.problemF \\
\kw{input:} \\
\> \( \omega = [\mathbf{a}, \mathbf{b}] \): an interval; \\
\> \( n \): an integer; \\
\kw{begin} \\
\> \( \mathbf{c} :\approx (\mathbf{a} + \mathbf{b}) / 2 \); \\
\> \( Q \).enqueue (\( [\mathbf{a},\mathbf{c}] \), \( n \)); \\
\> \( Q \).enqueue (\( [\mathbf{c},\mathbf{b}] \), \( n \)); \\
\kw{end.}
\end{tabbing}
\end{algorithm}

The problem reported in Algorithm~\ref{alg:interval} by a call to the function \emph{controller.innerEmpty}, however, is of different nature, and directly related to the situation described in \textbf{(P2b)} in Section \ref{sec:assigning}.
The interval must be then abandoned (moved to \( \mathcal{P}^- \)).
We do not provide pseudocode for this algorithm, because it is trivial.

% ====================================================

\subsection{Subdivisions of parameter intervals}
\label{sec:chop}

In this subsection we introduce an algorithm for subdividing a parameter interval \( \omega \)
when the numerical computations indicate that \( \omega_n \) might intersect the critical neighbourhood \( \Delta \).
The purpose of this subdivision is to cut out the part of \( \omega \) that falls onto \( \Delta \),
and to leave as much as possible from the interval \( \omega \) in the form of one or two subintervals of \( \omega \) that can be iterated further.

We begin by introducing Algorithm~\ref{alg:bisection} that uses the idea of the bisection method explained in Section~\ref{sec:avoid} to find a possibly small (or large) value of the parameter \( a \) such that \( c_n (a) \) is proved numerically to be below (or above, depending on which one is requested, the parameter called \textit{below} is used to make the choice) a certain ``border'' value \( \mathbf{v} \). The course of action of the algorithm depends on whether \( c_n \) is increasing or decreasing, and this information is passed to the algorithm in the parameter called \textit{increasing}. The number of bisection steps to conduct is given by the parameter \( s > 0 \). The features of the algorithm are precisely stated in Proposition~\ref{prop:bisection} below.

% Nonindent is necessary here for hyperref to produce a correct hyper link, see explanation here: https://tex.stackexchange.com/questions/14000/enumerate-after-label
\begin{algorithm}\noindent\label{alg:bisection}
\rm
\begin{tabbing}
\hspace{.5 cm}\= \hspace{.5 cm}\= \hspace{.5 cm}\= \hspace{.5 cm}\= %
\hspace{.5 cm}\= \hspace{.5 cm}\= \hspace{.5 cm}\= \hspace{.5 cm}\= \kill
\kw{function} bisection \\
\kw{input:} \\
\> \( \omega = [\mathbf{a}, \mathbf{b}] \): an interval; \\
\> \( \mathbf{v} \): real number; \\
\> \( n, s \): positive integers; \\
\> \( \textit{increasing}, \textit{below} \): boolean values (true or false); \\
\kw{begin} \\
\> \kw{repeat} \( s \) \kw{times}: \\
\>\> \( \mathbf{m} :\approx (\mathbf{a} + \mathbf{b}) / 2 \); \\
\>\> \( \mathbf{c}^- :\leq c_n (\mathbf{m}) \); \\
\>\> \( \mathbf{c}^+ :\geq c_n (\mathbf{m}) \); \\
\>\> \kw{if} \( \textit{increasing} \) \kw{and} \( \textit{below} \) \kw{then} \\
\>\>\> \kw{if} \( \mathbf{c}^+ \leq \mathbf{v} \) \kw{then} \( \mathbf{a} := \mathbf{m} \); \kw{else} \( \mathbf{b} := \mathbf{m} \); \\
\>\>\> \( \mathbf{p} := \mathbf{a} \); \\
\>\> \kw{if} \( \textit{increasing} \) \kw{and} \kw{not} \( \textit{below} \) \kw{then} \\
\>\>\> \kw{if} \( \mathbf{c}^- \geq \mathbf{v} \) \kw{then} \( \mathbf{b} := \mathbf{m} \); \kw{else} \( \mathbf{a} := \mathbf{m} \); \\
\>\>\> \( \mathbf{p} := \mathbf{b} \); \\
\>\> \kw{if} \kw{not} \( \textit{increasing} \) \kw{and} \( \textit{below} \) \kw{then} \\
\>\>\> \kw{if} \( \mathbf{c}^+ \leq \mathbf{v} \) \kw{then} \( \mathbf{b} := \mathbf{m} \); \kw{else} \( \mathbf{a} := \mathbf{m} \); \\
\>\>\> \( \mathbf{p} := \mathbf{b} \); \\
\>\> \kw{if} \kw{not} \( \textit{increasing} \) \kw{and} \kw{not} \( \textit{below} \) \kw{then} \\
\>\>\> \kw{if} \( \mathbf{c}^- \geq \mathbf{v} \) \kw{then} \( \mathbf{a} := \mathbf{m} \); \kw{else} \( \mathbf{b} := \mathbf{m} \); \\
\>\>\> \( \mathbf{p} := \mathbf{a} \); \\
\> \kw{return} \( \mathbf{p} \); \\
\kw{end.}
\end{tabbing}
\end{algorithm}

\begin{proposition}\label{prop:bisection}
Let \( \omega = [\mathbf{a}, \mathbf{b}] \subset [1,\infty) \) be a compact interval.
Let \( n > 0 \) be an integer such that \( c_n \) is monotone on \( \omega \).
Let the constant called ``\textit{increasing}'' have the value ``true'' if and only if \( c_n \) is increasing.
Let \( s > 0 \) be an integer.
Let \( \mathbf{v} \in \underline{\omega_n} \).
Let \( \mathbf{p}^- \) be the number returned by Algorithm~\ref{alg:bisection} with the parameter ``\textit{below}'' set to ``true'',
and let \( \mathbf{p}^+ \) be the number returned by Algorithm~\ref{alg:bisection} with the parameter ``\textit{below}'' set to ``false''.

Then
\begin{align}
\label{eq:bisection}
& c_n ( \mathbf{p}^- ) \leq \mathbf{v} \leq c_n ( \mathbf{p}^+ ),
\end{align}
and the same holds true for the numerically computed bounds for \( c_n ( \mathbf{p}^- ) \) and \( c_n ( \mathbf{p}^+) \).
\end{proposition}

\proof
Consider the case in which \( c_n \) is increasing (the case of \( c_n \) decreasing is analogous),
and thus assume \textit{increasing} is set to ``true.''
We shall prove that \( c_n ( \mathbf{p}^- ) \leq \mathbf{v} \) (the other part is analogous),
and thus assume \textit{below} is set to ``true.''
By the assumptions, \( c_n (\mathbf{a}) \leq \mathbf{v} \leq c_n (\mathbf{b}) \),
so \( \mathbf{a} \) is a good initial guess for \( \mathbf{p}^- \), but we are going to get a tighter bound.
In the loop repeated \( s \) times, the approximate midpoint \( \mathbf{m} \) of the interval \( [\mathbf{a},\mathbf{b}] \) is computed.
Then a lower bound \( \mathbf{c}^- \) and an upper bound \( \mathbf{c}^+ \) for the value of \( c_n (\mathbf{m}) \) are computed, and compared with \( \mathbf{v} \).
If it was proved numerically that \( c_n (\mathbf{m}) \leq \mathbf{v} \)
then the interval \( [\mathbf{a},\mathbf{b}] \) is replaced with \( [\mathbf{m},\mathbf{b}] \),
and \( \mathbf{m} \) becomes a new candidate for \( \mathbf{p}^- \).
Otherwise, \( \mathbf{a} \) remains a candidate for \( \mathbf{p}^- \), but we tighten the interval \( [\mathbf{a},\mathbf{b}] \) by replacing it with \( [\mathbf{a},\mathbf{m}] \). After this step, \( c_n (\mathbf{p}^-) \leq \mathbf{v} \), and the same holds true after the number of \( s \) steps.

The fact that the same inequalities hold true for the numerically computed bounds follows immediately from the fact that precisely these bounds are computed in Algorithm~\ref{alg:bisection} and the corresponding inequalities verified to obtain \eqref{eq:bisection}. We note, however, that the numerical method for computing these bounds must be identical each time, or otherwise this final conclusion may not hold true, due to rounding.
\qed

\begin{remark}\label{rem:bisection}
Since the interval \( [\mathbf{a},\mathbf{b}] \) is halved \( s \) times in Algorithm~\ref{alg:bisection}, the precision with which \( \mathbf{p}^- \) and \( \mathbf{p}^+ \) are estimated corresponds to \( 2^{-s} \) of the initial size of the interval \( \omega \).
\end{remark}

Next, we are going to introduce Algorithm~\ref{alg:hitDelta} that uses Algorithm~\ref{alg:bisection} to chop the interval \( \omega \) into three pieces: \( \omega^1 \), \( \omega^2 \), and \( \omega^3 \), with mutually disjoint interiors, such that \( c_n (\omega^i) \cap \Delta = \emptyset \) for \( i = 1, 2 \) unless \( \omega^i \) is degenerate (a singleton). Conceptually, if \( c_n (\omega) \) intersects \( \Delta \) then we cut out the part \( \omega^3 \) that hits \( \Delta \) from the middle of \( \omega \), so that we can continue iterating the remaining two subintervals of \( \omega \). The two computed subintervals \( \omega^1 \) and \( \omega^2 \) are added to the queue, unless they are too small, which is defined in terms of a certain fraction of the length of the interval \( \omega \), and the controller is notified about the interval \( \omega^3 \) excluded from further computations.

% Nonindent is necessary here for hyperref to produce a correct hyper link, see explanation here: https://tex.stackexchange.com/questions/14000/enumerate-after-label
\begin{algorithm}\noindent\label{alg:hitDelta}
\rm
\begin{tabbing}
\hspace{.5 cm}\= \hspace{.5 cm}\= \hspace{.5 cm}\= \hspace{.5 cm}\= %
\hspace{.5 cm}\= \hspace{.5 cm}\= \hspace{.5 cm}\= \hspace{.5 cm}\= \kill
\kw{function} omega\_hit\_Delta \\
\kw{input:} \\
\> \( \omega = [\mathbf{a}, \mathbf{b}] \): an interval; \\
\> \( n \): an integer; \\
%\> \( \overline{\omega_n} = [\mathbf{u}^-, \mathbf{u}^+] \): an interval; \\
\> \( \underline{\omega_n} = [\mathbf{u}^-, \mathbf{u}^+] \): an interval; \\
\> \textit{increasing}: a boolean value (true or false); \\
\kw{begin} \\
\> \( [\mathbf{a}',\mathbf{b}'] := [\mathbf{a},\mathbf{b}] \); \\
\> let \( s > 0 \) be the number of bisection steps recommended by the controller; \\
\> \kw{if} \( \mathbf{u}^- < -\delta \) \kw{then} \\
\>\> \kw{if} \textit{increasing} \kw{then} \\
\>\>\> \( \mathbf{a}' \) := bisection ( \( \omega \), \( \mathbf{v} = -\delta \), \( n \), \( s \), \textit{increasing} = true, \textit{below} = true); \\
\>\> \kw{else} \\
\>\>\> \( \mathbf{b}' \) := bisection ( \( \omega \), \( \mathbf{v} = -\delta \), \( n \), \( s \), \textit{increasing} = false, \textit{below} = true); \\
\> \kw{if} \( \mathbf{u}^+ > \delta \) \kw{then} \\
\>\> \kw{if} \textit{increasing} \kw{then} \\
\>\>\> \( \mathbf{b}' \) := bisection ( \( \omega \), \( \mathbf{v} = \delta \), \( n \), \( s \), \textit{increasing} = true, \textit{below} = false); \\
\>\> \kw{else} \\
\>\>\> \( \mathbf{a}' \) := bisection ( \( \omega \), \( \mathbf{v} = \delta \), \( n \), \( s \), \textit{increasing} = false, \textit{below} = false); \\
\> \kw{if} \( \mathbf{a} \neq \mathbf{a}' \) \kw{then} \( Q \).enqueue (\( [\mathbf{a},\mathbf{a}'] \), \( n \)); \\
\> \kw{if} \( \mathbf{b} \neq \mathbf{b}' \) \kw{then} \( Q \).enqueue (\( [\mathbf{b}',\mathbf{b}] \), \( n \)); \\
\> controller.notify\_excluded\_interval (\( [\mathbf{a}',\mathbf{b}'] \)); \\
\kw{end.}
\end{tabbing}
\end{algorithm}

\begin{proposition}\label{prop:disjoint}
Let \( \omega = [\mathbf{a}, \mathbf{b}] \subset [1,\infty) \) be a compact interval of parameters.
Let \( n > 0 \) be an integer.
Assume that \( \underline{{\omega}_{n}}\subseteq {\omega}_{n}\subseteq \overline{{\omega}_{n}} \),
and that \( \overline{{\omega}_{n}} \cap \Delta \neq \emptyset \).
Assume \( c_n \) is monotone on \( \omega \) and the value of the parameter ``\textit{increasing}'' is \emph{true} if and only if \( c_n \) is increasing.

Then all the intervals \( \omega^i \) added to the queue \( Q \)
by Algorithm~\ref{alg:hitDelta} applied to these objects
satisfy the following:
\begin{align}
& \omega^i \subset \omega. \\
& \overline{\omega^i_n} \cap \Delta = \emptyset,
\end{align}
where \( \overline{\omega^i_n} \) is the interval computed in Algorithm~\ref{alg:interval} for \( \omega^i \).
\end{proposition}

\proof
By construction, it is obvious that \( \omega^i \subset \omega \).
We are going to prove that the outer bound \( \overline{\omega^1_n} \) for \( c_n ([\mathbf{a},\mathbf{a}']) \) does not intersect \( \Delta \) if \( \mathbf{a} \neq \mathbf{a}' \) (the argument about \( [\mathbf{b}',\mathbf{b}] \) is analogous). By Proposition~\ref{prop:bisection}, if \( c_n \) is increasing on \( \omega \) and \( \mathbf{u}^- < -\delta \) then the bisection method provides \( \mathbf{a}' \) for which the numerically computed upper bound \( \mathbf{w} \) for \( c_n (\mathbf{a}') \) satisfies \(\mathbf{w} \leq -\delta \), and then indeed \( c_n ([\mathbf{a},\mathbf{a}']) \cap \Delta = \emptyset \), also as computed in the numerical version (with rounding).
\qed

% ====================================================

\subsection{Software}
\label{sec:software}

A software implementation of the algorithms introduced above  is publicly available at \cite{software}. The program is a command-line utility (to be launched in a text terminal, or at the command prompt), written in C++. It complies with the GNU C++ compiler (version 8.3.0, as of writing the paper). The GNU MPFR software library \cite{mpfr} is used for arithmetic operations on real numbers whenever high precision of the results and conrolled rounding are necessary. In particular, all real numbers provided in the input in the decimal form are rounded to the nearest representable numbers at the target precision by an MPFR function.

% In order to save the less computer skilled the hassle of downloading and compiling the source code ...:-)
We additionally provide a web interface at \cite{software} that makes it possible to run the program and see the results directly from the web browser. One fills out a table with the arguments to be passed to the program, hits the button to submit the form, and obtains the output produced by the program directly in the web page. The web interface allows the user to specify the parameter interval of interest and set several of the parameters involved in the computations, such as \( \delta, N_{0} \) and other parameters discussed in Section \ref{sec:assigning}, and possibly some other parameters not documented here (related, for example to the form in which the output is presented). For further details, we refer the reader to~\cite{software}.

% ====================================================

\section*{Acknowledgments}
\label{sec:ack}

A.G. and C.E.K. are grateful to ICTP, where part of this research was carried out, for its generous hospitality.

% ====================================================

\section*{Data availability statement}
\label{sec:data}

The raw data generated by our software, which constitutes a basis for the figures and for the analysis conducted in Section~\ref{sec:results}, is available in~\cite{data20}.

% ====================================================

% ====================================================


\begin{thebibliography}{99}

\bibitem{ArbMat04} 
A. Alexander and C. Matheus, ``Decidability of Chaos for Some Families of Dynamical Systems,'' \emph{Foundations of Computational Mathematics}  4 (3), 269--275 (2004).

\bibitem{AviLyudMel03} 
Avila, Artur, Mikhail Lyubich, and Welington de Melo, ``Regular or Stochastic Dynamics in Real Analytic Families of Unimodal Maps.'' \emph{Inventiones Mathematicae} 154 (3), 451--550 (2003)

\bibitem{BenCar85} M. Benedicks and L. Carleson,  ``On iterations of \( 1 - ax^{2}\) on (-1, 1)'' \newblock{\em Annals of Math.} 122: 1-25 (1985)

\bibitem{BenCar91}
M.~Benedicks and L.~Carleson, ``The dynamics of the H{\'e}non map,''
\newblock {\em Annals of Math.}, 133:73--169 (1991).

\bibitem{DKLMOP08}
S. Day, H. Kokubu, S. Luzzatto, K. Mischaikow, H. Oka and P. Pilarczyk,
\newblock{``Quantitative hyperbolicity estimates in one-dimensional dynamics,''}\newblock{\em Nonlinearity} 21: 1967--87 (2008).

\bibitem{Gal17}
Z.~Galias, ``Systematic search for wide periodic windows and bounds for the set of regular parameters for the quadratic map.'' \emph{Chaos} 27, 053106 (2017).

\bibitem{GolLuzPil16}
A. Golmakani, S. Luzzatto and P. Pilarczyk, \newblock{``Uniform expansivity outside a critical neighborhood in the quadratic family,''} \newblock{\em Exp. Math.} 25(2), 116--124 (2016).

\bibitem{GraSwi97} J. Graczyk and G. Swiatek,  \newblock{``Generic hyperbolicity in the logistic family,''} \newblock {\em Annals of Math.} 146:1-52 (1997).

\bibitem{mpfr}
% https://www.mpfr.org/algo.html
L. Fousse, G. Hanrot, V. Lef\`{e}vre, P. P\'{e}lissier, and P. Zimmermann.
``MPFR: A multiple-precision binary floating-point library with correct rounding,''
\emph{ACM Transactions on Mathematical Software,} 33(2), (2007).

\bibitem{ieee754}
IEEE 754--2019 - IEEE Standard for Floating-Point Arithmetic.
https://standards.ieee.org/content/ieee-standards/en/standard/754-2019.html.
Accessed on March 11, 2020.

\bibitem{Jak81} M. V. Jakobson \newblock{``Absolutely continuous invariant measures for one parameter families of one dimensional maps,''}\newblock{\em Commun. Math. Phys.} 81: 39--88 (1981).

\bibitem{Jak01} M. Jakobson \newblock{``Piecewise smooth maps with absolutely continuous invariant measures and uniformly scaled Markov partitions,''}\newblock{\em Smooth ergodic theory and its applications} (Seattle, WA, 1999), 825--881, Proc. Sympos. Pure Math., 69, Amer. Math. Soc., Providence, RI, 2001.

\bibitem{KinChe1991}
D. Kincaid and W. Cheney, ``Numerical analysis. Mathematics of scientific computing,''
Brooks/Cole Publishing Company, Pacific Grove, CA (1991).

\bibitem{LuzTuc99}
S.~Luzzatto and W.~Tucker,
\newblock {``Non-uniformly expanding dynamics in maps with criticalities and
  singularities,''} \newblock {\em Publ.  Math.  IHES}, 89: 179--226 (1999).

\bibitem{LuzVia00}
S.~Luzzatto and M.~Viana,
\newblock {``Positive Lyapunov exponents for Lorenz-like maps with
  criticalities,''} \newblock {\em Ast{\'e}risque}, 261: 201--237 (2000).

\bibitem{LuzTak06}
S. Luzzatto and H. Takahashi,
\textit{``Computable starting conditions for the existence of nonuniform hyperbolicity in one-dimensional maps,''} \emph{Nonlinearity} 19,  1657--1695 (2006).

\bibitem{Lyu97} M. Lyubich, \newblock{``Dynamics of quadratic polynomials I,'' } \newblock{\em Acta Math.} 178: 185-247 (1997).

\bibitem{Lyu97a} M. Lyubich, \newblock{``Dynamics of quadratic polynomials, II'' } \newblock{\em Acta Math.} 178: 247-97 (1997).

\bibitem{Lyu02}
M. Lyubich 
\newblock{``Almost Every Real Quadratic Map Is Either Regular or Stochastic,''}
\newblock{\emph{Annals of Math.} 156 (1) (2002). }

\bibitem{Moore1966}
R.E. Moore.
\textit{Interval analysis}.
Prentice-Hall, Inc., Englewood Cliffs, N.J. (1966).

\bibitem{PacRovVia98} M.J. Pacifico, A. Rovella, and M. Viana, 
\newblock{``Infinite-modal maps with global chaotic behavior,''}
\newblock{\emph{Annals of Math.} 148 (2) 441--484 ,(1998). } 
   
\bibitem{software}
P. Pilarczyk. \textit{Quadratic map software}.
\href{http://www.pawelpilarczyk.com/quadr/}{http://www.pawelpilarczyk.com/quadr/} (accessed on June 11, 2020).

\bibitem{data20}
P. Pilarczyk,
\newblock{``Stochastic intervals for the family of quadratic maps,''}
Gda\'{n}sk University of Technology (2020),
\href{https://dx.doi.org/10.34808/szfn-gv40}{doi: 10.34808/szfn-gv40}.

\bibitem{TucWil09} 
W. Tucker and D. Wilczak 
``A rigorous lower bound for the stability regions of the quadratic map,''
\emph{Physica D} 1923--1936 (2009).

\bibitem{WT2011}
W. Tucker, ``Validated numerics: a short introduction to rigorous computations'', Princeton University Press (2011).

\end{thebibliography}
\end{document}